\RequirePackage[reqno]{amsmath} 
\documentclass[11pt,reqno]{amsart}

\usepackage[foot]{amsaddr}
\usepackage{a4wide} 
\usepackage{mypub}

\usepackage{subfiles}

\begin{document}

\title[Continuity of the flow map]{Continuity of the flow map for
  symmetric hyperbolic systems and its application to the
  Euler--Poisson system}

  \author[U. 
  Brauer]{Uwe Brauer}

  \address{%
    Uwe Brauer Departamento de Matemática Aplicada\\ Universidad
    Complutense Madrid 28040 Madrid, Spain} \email{oub@mat.ucm.es}

  \thanks{U.~B. 
    gratefully acknowledges support from Grant MTM2016-75465 by
    MINECO, Spain and UCM-GR17-920894.} 

\author[L. 
  Karp]{Lavi Karp}

  \address{%
    Lavi Karp
    Department of Mathematics\\ ORT Braude College\\
    P.O. 
    Box 78, 21982 Karmiel\\ Israel}

  \email{karp@braude.ac.il}




\begin{abstract}
  We show the continuity of the flow map for quasilinear symmetric
  hyperbolic systems with general right--hand sides in different
  functional setting, including weighted Sobolev spaces $H_{s,\delta}$. 
  An essential tool to achieve the continuity of the flow map is a 
  new type of energy estimate, which we all it a low
  regularity energy estimate. 
  We then apply these results to the Euler--Poisson system which
  describes various systems of physical interests. 
\end{abstract}
\maketitle{}

\def\biblio{}

\section{Introduction}
\label{sec:introduction}

The purpose of our work is to prove the continuity of the flow map for
quasilinear symmetric hyperbolic systems in the topology of either the
ordinary Sobolev spaces~$H^s$, or the weighted Sobolev spaces
$H_{s,\delta}$, and under the assumption that lower order terms have limited
regularity. 
Once that is proven, we then apply this result to the
Euler--Poisson--Makino system. 
The Euler--Poisson--Makino system is a modification of the Euler--Poisson system
in which the density is replaced by a variable, which we denote as the
Makino variable from now on. 
This variable is a nonlinear function of the density. 
This variable change allows us to include situations where the density
can be zero. 
The Euler--Poisson--Makino system consists of quasilinear symmetric
hyperbolic evolution equations coupled to an elliptic equation.
The existence and the uniqueness of solutions in this setting have been
proved already by Makino
\cite{makino_86} in the $H^s$ spaces, and recently by
\cite{BK9} in the weighted Sobolev spaces $H_{s,\delta}$. 
Thus our result about the continuity of the flow map shows that these
systems are well-posed in the sense of Hadamard.

The general  symmetric hyperbolic system we have in mind has the
following form:
\begin{equation}
  \label{eq:initial-value:1}
  \begin{cases}
    A^0(U)\partial_t U+ \sum_{a=1}^{d} A^a(U)\partial_a U  =G(U), \\
    U(0,x)                                 =u_0(x),
  \end{cases}
\end{equation}
where $A^0,A^1,\ldots,A^d$ are $N\times N$ smooth symmetric matrices,
$A^0$ is positive definite and $G:\setR^N\to \setR^d$ is a nonlinear function.
It is well known that if $u_0$ belongs to the Sobolev space $H^s$ and
$s>\frac{d}{2}+1$, then there exists a positive $T$ and a unique
solution $U$ to \eqref{eq:initial-value:1} such that
\begin{equation}
  \label{eq:existence:1}
  U\in C^0([0,T];H^s)\cap C^1([0,T];H^{s-1}).
\end{equation}

This result was first proved by Marsden and Fischer \cite{FMA}, and
Kato \cite{KATO}. 

The goal of this paper is to investigate the continuous dependence on
the initial data, or equivalently the continuity of the flow map. 

\begin{defn}[Continuity of the flow map]
  We say that the flow map is {\bf continuous}, if for any $u_0\in H^s$, there
  exists a neighborhood $B\subset H^s$ of $u_0$ such that for every
  $u\in B$ the map $u\mapsto U$ from $B$ to $C^0([0,T];H^s)$ is continuous,
  where $U$ denotes  the solution to \eqref{eq:initial-value:1} with initial
  data $u$.
\end{defn}

Equivalently, let $\{u_n\}\subset B$ and $U^n$ be the corresponding solution to 
  \eqref{eq:initial-value:1} with initial data $u_n$, and $U^0$ the solution with 
initial data $u_0$. Then the flow map is 
continuous, if $u_n\to u_0$ in $H^s$ implies that $U^n\to U^0$ in $C^0([0,T];H^s)$.

Kato developed the abstract theory of nonlinear semigroups, which
also can be applied for many other types of evolution equations
\cite{kato75:_quasi} and which provides a tool to prove the continuity of the flow map. 
Subsequently, simpler methods that require elementary existence theory
of linear equations were implemented by various authors
\cite{majda84:_compr_fluid_flow_system_conser}, \cite{taylor97c}, \cite{Bahouri_2011}. 
However, those methods did not consider the continuity of the flow map
with one exception which we will discuss below.

The continuity of the flow map in a lower norm is relatively easy to
obtain (see, e.g., \cite[Thm.~6.6.5 p.~239]{rauch12:_hyper}); however,
it is unsatisfactory as a final result. 
The essential difficulty is to prove the continuity with respect to
the $H^s$ Norm. 
The main tool of our approach is a new type of energy estimate, namely,
one with lower regularity.
This new estimate, Lemma \ref{lem:sec1-outline:2}, concerns the
linearization of equation \eqref{eq:initial-value:1},
\begin{equation}
\label{eq:linear}
 A^0(t,x)\partial_t U+\sum_{a=1}^d A^a(t,x)\partial_a U=F(x,t).
\end{equation}

In the classical energy estimate all the components of
\eqref{eq:linear} are considered in the same Sobolev space $H^s$ (see
Lemma \ref{lem:sec1-math-prel:4}), while in the low regularity energy
estimate the solutions are in $H^{s-1}$, but the coefficients 
$A^0,A^1,\ldots,A^d$ belong to~$H^s$, and the range of the regularity
index $s$ is the same as in the standard one, namely,
$s>\frac{d}{2}+1$. 
We obtain this estimate by applying, in a slightly different way,  the
Kato--Ponce commutator estimate.

So one goal of this paper is to provide an alternative proof of Kato's
continuity result, a proof that is based on similar techniques used
for the proof of the existence and uniqueness theorems by elementary
methods.
We require a slightly different assumption than Kato, namely, we
demand that $DG$ is a Lipschitz function with respect to the $H^{s-1}$
norm (see Theorem \ref{thm:1}). 
This condition is almost equivalent to Lipschitz continuity in the
$H^s$ norm (see Remark \ref{rem:Lipschitz}).

In the book \cite{Bahouri_2011} the authors proved, besides existence and {unique-\linebreak ness}~\eqref{eq:existence:1}, the continuity of the flow map with zero
right-hand side of \eqref{eq:initial-value:1}.
We follow their ideas, which consist in analyzing the derivatives of
the solution, a suitable splitting of these derivatives which leads to 
two different systems, and then proving the continuous dependence,
using these two systems, in the $H^{s-1}$ norm (see the proof outline
following Remark~\ref{rem:Katos-example}). 
However, this approach relies on a certain type of energy estimate
which is not included in their work. 
That is why we have established this new type of estimate, that we
call a \textbf{low regularity energy estimate} (Lemma
\ref{lem:sec1-outline:2}) and which allows us to close the gap of the
missing estimate and therefore to complete the proof of the continuity
of the flow map for the systems we consider.

Moreover, since the authors only considered homogeneous systems
without source terms, we will, therefore, generalize their method to
include these terms, which play an important role in the coupled
Euler--Poisson system.

Another advantage of our approach is its extension of the continuity
of the flow map to the weighted Sobolev spaces $H_{s,\delta}$, which are
defined at the beginning of Section \ref{sec:h_s-delta-spaces} (see
Definition \ref{def:weighted:3}). 
Here $s$ denotes the regularity index as in the $H^s$ spaces, and
$\delta$ is the weight's index. 
These weighted spaces have the property that the weights vary with
respect to the order of the derivatives. 
Nirenberg and
Walker~\cite{nirenberg73:_null_spaces_ellip_differ_operat_r}, and
independently Cantor~\cite{cantor75:_spaces_funct_condit_r},
introduced these spaces for an integer order regularity index $s$. 
Triebel extended them to fractional order and proved basic properties
such as duality, interpolation, and density of smooth
functions~\cite{triebel76:_spaces_kudrj2}.

Traditionally, symmetric hyperbolic systems have been treated in
Sobolev spaces $H^s$ because their norm allows in a convenient way to
obtain energy estimates. 
But there are situations in which these spaces are too restrictive
\cite{KATO, majda84:_compr_fluid_flow_system_conser},
for example, in the case of the Euler--Poisson system, if the density
has no compact support, or the Einstein equations in asymptotically
flat spaces. 
For these types of problems, the weighted Sobolev spaces $H_{s,\delta}$ are
useful.
The advantage of these spaces is that we can control the falling-off/growth rate near infinity. 
These spaces have been successfully used to treat many elliptical and
geometrical problems (see, for example, \cite{choquet--bruhat81:_ellip_system_h_spaces_manif_euclid_infin} 
and \cite{Lee_Parker_87}).

Recently the authors established existence and uniqueness theorems for
the system \eqref{eq:initial-value:1} in the $H_{s,\delta}$ spaces, \cite{BK9}. 
More precisely, if $u_0$ belongs to $H_{s,\delta}$, $s>\frac{d}{2}+1$, and
$\delta\geq -\frac{d}{2}$, then there is a positive $T$ and a unique solution
to \eqref{eq:initial-value:1} such that
\begin{equation*}
 U\in C^0([0,T]; H_{s,\delta})\cap C^1([0,T]; H_{s-1,\delta+1}). 
\end{equation*}
In the present paper we prove the continuity of the flow map in the
topology of the $H_{s,\delta}$ norm under the condition that $DG$ is
Lipschitz in the $H_{s-1,\delta}$ norm. 
Hence we conclude that the initial value problem for first order
symmetric hyperbolic systems in the $H_{s,\delta} $ spaces is well-posed in
the sense of Hadamard. 
This result we apply to the Euler--Poisson--Makino system for densities that
are not compactly supported.

The problem of continuous dependence on the data in the strong
topology can be solved rather easily for elliptic and parabolic type
partial differential equations. 
But it is considerably more complicated for initial value problems for
fluids and hyperbolic equations as  pointed out by Kato and Lai
\cite{Kato_Lai_84}.
For the incompressible Euler equations the continuity of the flow map
was proven in the celebrated paper by Kato and Ponce
\cite{kato88:_commut_euler_navier_stokes}. 
Beir\~{a}o da
Veiga~\cite{beirao92:_data} proved the continuity of the flow map for the compressible Euler
equations in a setting of an initial boundary value problem.
Speck \cite{Speck_09} investigated the well-posed problem in the Hadamard sense for the simplified
Euler--Nordstr\"om system. 
In order to do so, he wrote the system as a first order hyperbolic
system. 
He did not rely on the theory of symmetric hyperbolic systems;
instead, he used the theory of energy currents which was introduced by
Christodoulou
\cite{christodoulou00A}.

For certain types of evolutionary equations, for example, the
Korteweg--de Vries (KdV) and nonlinear Schr\"odinger equations, the
flow map is Lipschitz continuous
(see~\cite{tzvetkov06:_ill} and the reference therein). 
But for the most simple quasilinear symmetric hyperbolic equation,
namely, the Burgers equation
\begin{equation*}
 \partial_t U+U\partial_x U=0,
\end{equation*}
Kato~\cite{KATO} demonstrated that the flow map cannot be of H\"older continuous type.
This example indicates that the continuity of the flow map for a
symmetric hyperbolic system is a delicate issue.

The paper proceeds as follows: In
Section~\ref{sec:continuity-flow-map-4} we recall some basic
notations; this is followed by Subsection~\ref{sec:useful-lemmas}
where we prove the low regularity energy estimate in the $H^s$ norm. 
This estimate allows us to prove the continuous dependence of the
solutions in the coefficients. 
Based on this result we can then show in Subsection~\ref{sec:continuity-flow-map} the continuity of the flow map for
solutions of quasilinear equations in the $H^s$ spaces. 
The next section, Section \ref{sec:continuity-flow-map-1}, has
basically the same structure but using weighted Sobolev spaces
$H_{s,\delta}$. 
In Section \ref{sec:euler-poisson-makino} we apply our results to
situations of physical interest. 
In Subsection \ref{sec:case-comp-supp} we apply it to a model in which
the density has compact support, that is followed, in Subsection
\ref{sec:case-density-which}, by a model whose density falls off at
infinity in an appropriate manner. 
The last example, in Subsection \ref{sec:cosmological-context},
concerns a cosmological situation in which the density extends all
over the space. 

Finally, in the appendix, we provide some tools which are useful for our
purpose and which we present, for the convenience of the reader, at the
end of the paper.

\biblio

\section{Continuity of the flow map in $H^s$}
\label{sec:continuity-flow-map-4}

\subsection{Notation and symbols.}
\label{sec:notations-symbols}

Let
\begin{equation}
  \label{eq:notation:1}
  \Lambda^s[u]=\mathcal{F}^{-1}(1+|\xi|^2)^{\frac{s}{2}}\mathcal{F}[u],
\end{equation}
where $\mathcal{F}[u](\xi)=\widehat{u}(\xi)$ is the Fourier transform. 
The $H^s$ norm is defined by
\begin{equation*}
  \| u\|_{H^s}= \|\Lambda^s[u]\|_{L^2}=\bigg(c_d\int_{\setR^d} 
    (1+|\xi |^2)^s|\widehat{u}(\xi)|^2d\xi\bigg)^{\frac{1}{2}}.
\end{equation*}
We will denote $H^s=H^s(\setR^d;\setR^N)$ the space of all functions
$u:\setR^d\to\setR^N$ with a finite $H^s$ norm. 
This is a Hilbert space with an inner product
\begin{equation}
  \label{eq:inner-product}
  \langle u,v\rangle_s=\langle \Lambda^s[u], 
  \Lambda^s[v]\rangle_{L^2}=\int_{\setR^d} (\Lambda^s[u](x)\cdot \Lambda^s[v](x)) dx,
\end{equation}
where $\cdot$ denotes the scalar product.

We will use the symbol $D=D_x$ to denote the derivative with respect to the
space variable.

\subsection{Energy estimates in the $H^s$ spaces.}
\label{sec:useful-lemmas}

In this section, we consider energy estimates for a linear symmetric
hyperbolic system
\begin{equation*}
  \label{eq:sec1-math-prel:3}
  \tag{LS}
\begin{cases}
           \partial_t U+ \sum_{a=1}^d A^{a}(t,x)\partial_a U =F(t,x)\\
           U(0,x)=u_0(x).
         \end{cases}
\end{equation*}
where $A^a$ are $N\times N $ symmetric matrices. The existence and the uniqueness
of solutions in the $H^s$ space are well-known (see, e.g., \cite{Bahouri_2011},
\cite{rauch12:_hyper}) and that  is why we will not  mention it in
the following every time where it seems necessary.

We will prove a new energy estimate with lower regularity for functions
\mbox{$U(t)=U(t,x)$} which are solutions to the linear systems
\eqref{eq:sec1-math-prel:3}. 
In order to emphasize the difference between the two types of
estimates, we start with the presentation of the traditional energy
estimate.

\begin{lem}[Standard energy estimate]
  \label{lem:sec1-math-prel:4}
  Let
   $s>\frac{d}{2}+1,$
   assume $A^a$, \mbox{$F\in L^\infty([0,T]; H^{s}) $} and
  $u_0\in H^s$. 
  Then a solution $U(t) $ to the linear system
  \eqref{eq:sec1-math-prel:3} that belongs to
  $ C^{0}( [0,T]; H^{s}) \cap C^{1}( [0,T]; H^{s-1})$
  satisfies the estimate
  \begin{equation*}
    \label{eq:sec2-lemmas:1}
    \Vert U(t) \Vert_{H^{s}}^2\leq e^{\int_0^t a_s(\tau)d\tau} 
    \bigg( 
      \Vert u_0\Vert_{H^{s}}^2+ \int_0^t \Vert F(\tau)\Vert_{H^{s}}^2 d\tau\bigg),
  \end{equation*}
  where
  $a_s(\tau):= C\sum_{a=1}^d\Vert    A^a(\tau)\Vert_{H^s}+1$. 
\end{lem}

We turn now to the low regularity  energy estimate.
\begin{lem}[Low regularity energy estimate]
  \label{lem:sec1-outline:2}
  Let $s\!>\!\frac{d}{2}\!+\!1$, $A^a\! \in\! L^\infty([0,T]; H^{s}) $,
  $F \in L^\infty([0,T]; H^{s-1}) $, and $u_0\in H^{s-1}$.
  Assume $U(t)\in L^{\infty}([0,T] ;H^{s-1})$ is a solution to
  the initial value problem \eqref{eq:sec1-math-prel:3}. Then for
  $t\in [0,T]$,
  \begin{equation}
    \label{eq:sec2-lemmas:15}
    \Vert U(t) \Vert_{H^{s-1}}^2\leq e^{\int_0^t a_s(\tau)d\tau} 
    \bigg( \Vert u_0 \Vert_{H^{s-1}}^2+ \int_0^t     \Vert 
        F(\tau,\cdot)\Vert_{H^{s-1}}^2 d\tau \bigg),
  \end{equation}
  where
  $a_s(\tau):=  C\sum_{a=1}^d\Vert
    A^a(\tau)\Vert_{H^s}+1$.  
\end{lem}

\begin{rem}
\label{rem:sec1-cont-flow-hs:1}
  Note that in Lemma \ref{lem:sec1-outline:2} we assume that the
  matrices $A^a$ have one more degree of differentiability than the
  other terms such as the initial data $u_0$, the right-hand side $F$
  and the solution $U(t)$. 
  On the other hand, in Lemma \ref{lem:sec1-math-prel:4} all those
  terms are in $H^s$.
\end{rem}

The proof of Lemma \ref{lem:sec1-outline:2}
contains a subtle point that occurs also in the proof of Lemma
\ref{lem:sec1-math-prel:4}. 
The common way to prove these types of energy estimates is to
differentiate $\|U(t)\|_{H^s}^2$ with respect to $t$, and that results in
the identity
\footnote{Strictly speaking, we have to justify this
  identity; however, for the sake of brevity, we leave it out but refer,
  for example, to the book of Temam \cite[Ch.~III, Lemma 1.1]{temam79:_navier_stokes} who has even treated a more general setting.}
\begin{equation}
 \label{eq:sec1:19}
 \frac{d}{dt}\|U(t)\|_{H^s}^2=2\langle U(t),\partial_t U(t)\rangle_s.
\end{equation} 
Assuming then that $U(t)$ is a solution to
\eqref{eq:sec1-math-prel:3}, we insert the equation into
\eqref{eq:sec1:19} and proceed by using integration by parts and
appropriate properties of Sobolev spaces; in the final step, we apply
the Gronwall inequality. 

However, the right-hand side of equation \eqref{eq:sec1:19},
$2\langle U(t),\partial_t U(t)\rangle_s$, is defined only if
$\partial_tU(t)\in H^s$, but on the other hand, if
$U(t)\in L^{\infty}( [0,T]; H^{s}) \cap C^{1}( [0,T];
  H^{s-1})$ is a solution to \eqref{eq:sec1-math-prel:3}, then
$\partial_t U(t)$ only belongs to $H^{s-1}$. 
Consequently, we can apply the identity \eqref{eq:sec1:19} only under
the assumption that the solutions belong to $H^{s+1}$. 
In order to justify this assumption, we have to perform two main steps.

In the first step, we assume that $u_0$ and $F$ have one degree of
regularity more, say $s+1$, and we establish the energy estimate in
the $H^{s}$ norm. 
In the second step, we approximate the terms $u_0$ and $F$ by sequences
in $H^{s+1}$ and show that the sequence of the corresponding solutions
converges weakly in the topology of $H^s$. 

Hence we start with the following proposition.
\begin{prop}
  \label{lem:sec1:2}
  Let $s>\frac{d}{2}+1$, $A^a \in L^\infty([0,T]; H^{s}) $,
  $F \in L^\infty([0,T]; H^{s}) $, $u_0\in H^{s}$. 
  Assume
  $U(t)\in C^{0}( [0,T]; H^{s}) \cap C^{1}( [0,T];
    H^{s-1})$ is a solution to the initial value problem
  \eqref{eq:sec1-math-prel:3}. Then for $t\in [0,T]$,
  \begin{equation}
    \label{eq:sec1:1}
    \Vert U(t) \Vert_{H^{s-1}}^2\leq e^{\int_0^t a_s(\tau)d\tau} 
    \bigg( \Vert u_0 \Vert_{H^{s-1}}^2+ \int_0^t \Vert 
        F(\tau,\cdot)\Vert_{H^{s-1}}^2 d\tau\bigg),
  \end{equation}
  where
  $a_s(\tau):= C \sum_{a=1}^d\Vert
    A^a(\tau)\Vert_{H^s}+1$.  
\end{prop}
\begin{proof}
  The proof follows the standard techniques as described above. 
  The essential new ingredient is a modification of the common
  application of commutator
  $[\Lambda^{s-1}A^a,A^a\Lambda^{s-1}](\partial_a U)$ and will be explained
  below. 
  Since $\partial_tU, \partial_a U\in H^{s-1}$, we can apply identity \eqref{eq:sec1:19}
  and conclude that
  \begin{equation}
    \label{eq:sec1:2}
    \begin{split}
       \frac{1}{2} \dfrac{d }{dt} \|U(t)\|_{H^{s-1}}^2 &=
      \langle U(t),\partial_t U(t)\rangle_{{s-1}}\\ 
      &= - \sum_{a=1}^d\langle
        U(t),A^a(t,\cdot)\partial_a U(t)\rangle_{{s-1}}+ \langle U(t),F(t,\cdot)\rangle_{s-1}.
    \end{split}
  \end{equation}
  We start with the last (the low order) term in \eqref{eq:sec1:2},
  which we estimate easily by the Cauchy--Schwarz inequality,
  \begin{equation}
    \label{eq:sec1:3}
    \begin{split}
      \Big|\langle U(t),F(t,\cdot)\rangle_{{s-1}}\Big| & \leq
      \|U(t)\|_{H^{s-1}}\|F(t,\cdot)\|_{H^{s-1}} \\
       & \leq \frac{1}{2}(\|U(t)\|_{H^{s-1}}^2+
        \|F(t,\cdot)\|_{H^{ s-1 } } ^2).
    \end{split}
  \end{equation}
  
  \noindent
  The next step consists in estimating the first order term which
  contains  the matrices~$A^a$.  
  So suppose we start by considering the following commutator:
  \begin{equation}
    \label{eq:sec1:4}
    \begin{split}
      \langle U,A^a\partial_a U\rangle_{{s-1}} & = \langle
        \Lambda^{s-1}[U],\Lambda^{s-1}[A^a\partial_a U]\rangle_{L^2} \\
      & = \langle \Lambda^{s-1}[U],\{ \Lambda^{s-1}[{A}^a\partial_a
            U]-{A}^a \Lambda^{s-1}[\partial_a
            U]\}\rangle_{L^2} \\ 
            &\phantom{=\ }+ \langle \Lambda^{s-1}[U],A^a
        \Lambda^{s-1}[\partial_a U]\rangle_{L^2}.
    \end{split}
  \end{equation}

  Then the last term on the right-hand side of \eqref{eq:sec1:4} is taken
  care of by integration by parts. 
  We turn now to the first term, use the Cauchy--Schwarz inequality,
  and then we need to estimate the $L^2$ norm of
  \begin{equation*}
    \{ \Lambda^{s-1}[{A}^a\partial_a U]-{A}^a
      \Lambda^{s-1}[\partial_a U]\}.
  \end{equation*}

  This term, however, presents some difficulties and that is why we
  will discuss this crucial issue in some detail. 
  The common way to estimate this commutator is by using the
  Kato--Ponce estimate \eqref{eq:kato-ponce} (see Taylor
  \cite[\S3.6]{Taylor91}).
  However, if we applied it directly, we would obtain
  \begin{equation}
\label{eq:sec1-cont-flow-hs:2}
    \begin{split}
       \|\{\Lambda^{s-1}[A^a\partial_a U]-A^a
          \Lambda^{s-1}[\partial_a &U]\}\|_{L^2} \\ \lesssim &
      \{\|D A^a\|_{L^\infty}\|\partial_a
          U\|_{H^{s-2}}+\|A^a\|_{H^{s-1}}\|\partial_a
          U\|_{L^\infty}\}.
    \end{split}.
  \end{equation}
  In order to obtain all the terms in the $H^{s-1}$ norm we use
  the Sobolev embedding, but this results in
  \begin{equation*}
    \|\partial_a U\|_{L^\infty}\lesssim \|\partial_a 
      U\|_{H^{s-2}}\lesssim \| U\|_{H^{s-1}},
  \end{equation*}
  which causes a loss of the regularity, since the Sobolev embedding theorem
  requires that $s-2>\frac{d}{2}$, while the assumption of the Lemma
  is $s-1>\frac{d}{2}$.

  The idea to overcome this apparent difficulty is to use a slightly
  different commutator. First, we write 
  \begin{equation*}
    A^a\partial_a U=\partial_a(A^a U)-\partial_a(A^a)U;
  \end{equation*}
  then \eqref{eq:sec1:4} becomes 
  \begin{equation}\label{eq:sec1-cont-flow-hs:3}
  \begin{aligned}
    \langle U, A^a\partial_aU\rangle_{s-1}    & 
=\langle\Lambda^{s-1}[U],\Lambda^{s-1} [A^a\partial_a U]\rangle_{L^2} \\
    &=\langle\Lambda^{s-1}[U], \Lambda^{s-1} [\partial_a ( A^aU ) 
]\rangle_{L^2} -\langle\Lambda^{s-1}[U], \Lambda^{s-1} [ \partial_a ( A^a
    )U ]\rangle_{L^2}.
  \end{aligned}
  \end{equation}
Now we form a commutator with respect to the first expression in
equation \eqref{eq:sec1-cont-flow-hs:3}, which results in 
\begin{equation}
\label{eq:sec1-cont-flow-hs:4}
\begin{split}
    \langle U,A^a&\partial_a U\rangle_{{s-1}}\\
    =\ &\langle\Lambda^{s-1}[U], \Lambda^{s-1} [\partial_a  ( A^aU 
]-A^a\Lambda^{s-1}[\partial_aU ]\rangle_{L^2}
    +  \langle\Lambda^{s-1}[U],  A^a\Lambda^{s-1}[\partial_a U] \rangle_{L^2}\\
     & -\langle\Lambda^{s-1}[U], \Lambda^{s-1} [ \partial_a ( A^a )U 
]\rangle_{L^2}.
   \end{split}
 \end{equation}

  The advantage of this expression is that we can now use the
  Kato--Ponce commutator estimate \eqref{eq:kato-ponce} with the
pseudo-differential operator $\Lambda^{s-1}\partial_a$, which belongs to class
  $OPS^s_{1,0}$, rather than $OPS^{s-1}_{1,0}$ as above. 
  So now we apply Kato--Ponce estimate \eqref{eq:kato-ponce} to the
  first term on the right-hand side of \eqref{eq:sec1-cont-flow-hs:4}, which
  results in 
  \begin{equation*}
    \label{eq:sec1:6}
       \|\{(\Lambda^{s-1}\partial_a)[A^a U]-A^a
          (\Lambda^{s-1}\partial_a)[ U]\}\|_{L^2}  \lesssim
       \{\|D A^a\|_{L^\infty}\|
          U\|_{H^{s-1}}+\|A^a\|_{H^{s}}\|
          U\|_{L^\infty}\}.
  \end{equation*}
  Now we apply the Sobolev embedding theorem to $\|U\|_{L^\infty}$, which allows
  us to conclude that $\|U\|_{L^\infty}\lesssim \|U\|_{H^{s-1}}$, for
  $s-1>\frac{d}{2}$. 
  If, on the other hand, we had used the commutator
  \eqref{eq:sec1-cont-flow-hs:2}, then  we would have to estimate the term
  $\|\partial_a U\|_{L^\infty}$ which would have
  contradicted the condition $s>\frac{d}{2}+1$. 
  That clearly shows the advantage of our approach.

  Estimating   $\|DA^a\|_{L^\infty}$ by the Sobolev embedding theorem
  results in
  \begin{equation}
    \label{eq:sec1:7}
    \|\{(\Lambda^{s-1}\partial_a)[A^a U]-A^a
        (\Lambda^{s-1}\partial_a)[ U]\}\|_{L^2}
    \lesssim \|A^a\|_{H^s}\|U\|_{H^{s-1}}.
  \end{equation}
  For the last term of \eqref{eq:sec1-cont-flow-hs:4} we use the multiplication
  property \eqref{eq:mult} which results in
  \begin{equation}
    \label{eq:sec1:8}
    \begin{split}
      \|\Lambda^{s-1}[(\partial_a A^a)U]\|_{L^2} 
      & =\|(\partial_a  A^a)U\|_{H^{s-1}}\lesssim \|\partial_a A^a\|_{H^{s-1}}
      \|U\|_{H^{s-1}} \\ 
      & \lesssim \| A^a\|_{H^{s}}
      \|U\|_{H^{s-1}}.
    \end{split}
  \end{equation}
  We turn now to the second term of \eqref{eq:sec1-cont-flow-hs:4} and apply
  integration by parts. Using the well-known identity
    \begin{equation*}
  \begin{split}
    \langle \Lambda^{s-1}[U],A^a\Lambda^{s-1}[\partial_aU]\rangle_{L^2} +
    \langle& \Lambda^{s-1}[\partial_a U], A^a\Lambda^{s-1}[U]\rangle_{L^2}
    \\&+\langle \Lambda^{s-1}[U],(\partial_a
    A^a)\Lambda^{s-1}[U]\rangle_{L^2}=0
  \end{split}.
\end{equation*}
    and the fact that  the matrices $A^a$ are symmetric, we can conclude that 
    \begin{equation}
      \label{eq:sec1:12}
\begin{split}
        |\langle \Lambda^{s-1}[U],&\,\,A^a[\Lambda^{s-1}\partial_a          U]\rangle_{L^2}|\\
       & = \frac{1}{2} \vert\langle 
\Lambda^{s-1}[U],(\partial_{a} A^a)\Lambda^{s-1}[U]\rangle_{L^2}  
\vert\\
        &\leq \|\partial_a A^a\|_{L^\infty}\|\Lambda^{s-1}[U]\|_{L^2}^2
        \lesssim \|A^a\|_{H^{s}}\|U\|_{H^{s-1}}^2.
      \end{split}
          \end{equation}
  Using inequalities \eqref{eq:sec1:3}, \eqref{eq:sec1:7},
  \eqref{eq:sec1:8} and \eqref{eq:sec1:12}, we end up with the energy
  estimate in differential form,
  \begin{equation}
    \label{eq:sec1:9}
    \frac{d}{dt}\|U(t)\|_{H^{s-1}}^2\lesssim 
    \bigg(\sum_{a=1}^d \|A^a(t,\cdot)\|_{H^{s-1}} +1\bigg) 
    \|U(t)\|_{H^{s-1}}^2 + \|F(t,\cdot)\|_{H^{s-1}}^2.
  \end{equation}
  We complete the proof by applying Gronwall's
  inequality~\eqref{eq:sec1-outline:53} to inequality \eqref{eq:sec1:9} and then
  inequality \eqref{eq:sec1:1} follows.
\end{proof}

It remains to prove Lemma \ref{lem:sec1-outline:2} for which we shall
use Proposition \ref{lem:sec1:2}, approximation, and weak convergence. 

\begin{proof}[Proof of Lemma \ref{lem:sec1-outline:2}] Since $H^s$ is dense
  in $H^{s-1}$, we take the following sequences:
  $F^k\in L^\infty([0,T];H^s)$ and $u_0^k\in H^s$ such that
  $$\|F^k(t,\cdot)-F(t,\cdot)\|_{H^{s-1}}\to 0\quad\text{and}\quad
  \|u_0^k-u_0\|_{H^{s-1}}\to 0, \quad\text{as } k\to\infty.$$
  Then for each $k$, there is a
  \mbox{$U^k\in C^0([0,T];H^s)\cap C^1([0,T];H^{s-1})$},
  which is a unique solution to the initial value problem
  \eqref{eq:sec1-math-prel:3} with right-hand side $F^k$ and initial
  data $u_0^k$ (see, e.g., 
  \cite[Theorem 4.4]{Bahouri_2011}). 
  Hence we can apply inequality \eqref{eq:sec1:1} of Proposition
  \ref{lem:sec1:2} to the sequence $\{U^k\}$ and obtain for
  sufficiently large $k$
\begin{equation}
    \label{eq:sec1:11}
    \begin{split}
      \Vert U^k(t) \Vert_{H^{s-1}}^2 &\leq e^{\int_0^t a_s(\tau)d\tau}
      \bigg( \Vert u_0^k \Vert_{H^{s-1}}^2+ \int_0^t\Vert
          F^k(\tau,\cdot)\Vert_{H^{s-1}}^2 d\tau \bigg)\\
     & \leq e^{\int_0^t a_s(\tau)d\tau} \bigg( \Vert u_0
        \Vert_{H^{s-1}}^2+1+ \int_0^t\Vert (
            F(\tau,\cdot)\Vert_{H^{s-1}}^2+1) d\tau \bigg).
    \end{split}
  \end{equation}
 Therefore for each $t\in[0,T]$, $\{U^k(t)\} $ is a bounded sequence in 
$H^{s-1}$, and hence converges weakly (up to the extraction of a  suitable 
subsequence) to $\widetilde{U}(t)\in H^{s-1}$. 
Note that $(U^k-U)$ satisfies the linear initial value problem 
\begin{equation*}
  \begin{cases}
 \partial_t(U^k-U) +  \sum_{a=1}^d A^a(t,x)\partial_a(U^k-U)=F^k(t,x)-F(t,x),\\
 (U^k-U)(0,x)=u_0^k(x)-u_0(x), 
\end{cases}
\end{equation*}
and then by the $L^2$ energy estimate we obtain 
\begin{equation*}
 \|U^k(t)-U(t)\|_{L^2}^2\leq e^{\int_0^t a_{\infty}(\tau)d\tau}\bigg(
 \|u_0^k-u_0\|_{L^2}^2+\int_0^t 
\|F^k(\tau,\cdot)-F(\tau,\cdot)\|_{L^2}^2d\tau\bigg),
\end{equation*}
where
$ a_\infty(\tau):=\sum_{a=1}^d\|\partial_a
A^a(\tau,\cdot)\|_{L^\infty}+1$.
Thus $U^k(t)\to U(t)$ in $L^2$ for each $t\in[0,T]$.
In conclusion, we have that $U^k(t)\to \widetilde{U}(t)$ weakly in
$H^{s-1}\subset L^2$ and $U^k(t)\to {U}(t)$ in the norm of $L^2$, hence by
Proposition \ref{prop:weak-conv} in  \ref{sec:mathematical-tools}, 
$\widetilde{U}(t)=U(t)$. 
So now by the weak limit and inequality \eqref{eq:sec1:11}, we obtain 
\begin{equation*}
\begin{split}
   \|U(t)\|_{H^{s-1}}^2&\leq 
\liminf_{k\to\infty}\|U^k(t)\|_{H^{s-1}}^2 
\\  & \leq \liminf_{k\to\infty}
e^{\int_0^t a_s(\tau)d\tau} 
\bigg( \Vert u_0^k \Vert_{H^{s-1}}^2+ \int_0^t\Vert 
F^k(\tau,\cdot)\Vert_{H^{s-1}}^2 d\tau \bigg)\\ 
&  =e^{\int_0^t a_s(\tau)d\tau} 
\bigg( \Vert u_0 \Vert_{H^{s-1}}^2+ \int_0^t\Vert 
F(\tau,\cdot)\Vert_{H^{s-1}}^2 d\tau \bigg)
\end{split},
\end{equation*} 
which completes the proof of the lemma.
\end{proof}

\subsection{Continuous dependence of the solution on its coefficients.}
\label{sec:cont-depend-solut}
The following lemma deals with the continuous dependence of the
solutions in~$H^{s-1}$ on the initial value problem
\eqref{eq:sec1-math-prel:3}. 
The continuous dependence is understood with respect to the matrices
$A^a$.
This lemma plays a central role in the proof of the main result.
\begin{lem}[Continuous dependence]
  \label{lem:sec1-outline:4}
  Let $s>\frac{d}{2}+1$,
  $ F\in L^{\infty}([0,T], H^{s-1} )$, $ u_0\in H^{s-1}$. 
  Let $A^a $ and $ \{A_n^a\}_{n=1}^\infty $ be symmetric matrices such that
  \begin{equation}
    \label{eq:sec1-outline:14}
    \Vert A^a(t,\cdot) \Vert_{H^s}  \leq C_0, \quad   
     \Vert A^a_n(t,\cdot) \Vert_{H^s}  \leq C_0, \quad t\in[0,T],   \; 
    \forall n\in\setN
  \end{equation}
  and
  \begin{equation}
    \label{eq:sec1-outline:54}
    \lim_{n\to\infty} \Vert A^a_n (t,\cdot)-A^a (t,\cdot)  \Vert_{L^\infty([0,T]; H^{s-1})} =0.
  \end{equation}
  Let $U(t), U^n(t)\in C^{0}( [0,T]; H^{s-1})$ be solutions to
  \eqref{eq:sec1-math-prel:3} with coefficients $A^a$ and~$A^a_n$,
  respectively, and with the same data $F$ and $u_0$. 
  Then
  \begin{equation*}
    \label{eq:sec1-cont-flow-hs:5}
    \lim_{n\to\infty}\sup_{0\leq t \leq T}\Vert U^n(t)- U(t) 
    \Vert_{H^{s-1}}=0.
  \end{equation*}
\end{lem}

\begin{proof}
  As in the proof of Lemma \ref{lem:sec1-outline:2}, we approximate
  $F$ and $u_0$ by sequences $F^k\in L^\infty([0,T]; H^s)$ and
  $u_0^k\in H^s$.
  Let $U^{n,k}$ be the solution of the initial value problem
  \eqref{eq:sec1-math-prel:3} with coefficients $A^a_n$, right-hand
  side $F^k$ and initial data $u_0^k$. 
  In addition, let $U^{k}$ be the solution of
  \eqref{eq:sec1-math-prel:3} with coefficients $A^a$, right-hand
  side $F^k$ and initial data $u_0^k$.
  We then make the following decomposition:
  \begin{equation}
    \label{eq:sec1:16}
    U^n(t)-U(t)=(U^n(t)-U^{n,k}(t))+(U^{n,k}(t)
      -U^k(t))+(U^k(t)-U(t)).
  \end{equation}
  We observe,  by Lemma \ref{lem:sec1-outline:2}, that
  \begin{equation}
    \label{eq:sec1:13}
    \begin{split}
      \Vert U^n(t)&-U^{n,k}(t) \Vert_{H^{s-1}}^2\\
      &\leq e^{\int_0^t
        a_{n,s}(\tau)d\tau} \bigg( \Vert u_0- u_0^k \Vert_{H^{s-1}}^2      
        + \int_0^t \Vert F(\tau,\cdot)-F^k(\tau,\cdot)\Vert_{H^{s-1}}^2
        d\tau\bigg),
    \end{split}
  \end{equation}
  where
  $a_{n,s}(\tau):= C\sum_{a=1}^d\Vert A^a_{n}(\tau)\Vert_{H^s}+1$. 
  Hence it follows by assumption \eqref{eq:sec1-outline:14} that
  $a_{n,s}\leq {C}$, which is independent of $n$. 
  That is why inequality \eqref{eq:sec1:13} implies that
  \begin{equation}
    \label{eq:sec:14}
    \lim_{k\to\infty}\sup_{0\leq t \leq T} \Vert U^n(t)-U^{n,k}(t) \Vert_{H^{s-1}}=0,
  \end{equation}
  uniformly in $n$.  Similarly
  \begin{equation}
    \label{eq:sec:15}
    \lim_{k\to\infty}\sup_{0\leq t \leq T}   \Vert U^k(t)-U(t) \Vert_{H^{s-1}}=0.
  \end{equation}
  Next, we note that $U^{n,k}-U^k$ satisfies the initial value problem
  \begin{equation*}
      \begin{cases}
         \partial_t(U^{n,k}-U^k) +\sum_{a=1}^d A^a_n(t,x)
        \partial_a(U^{n,k}-U^k)\\
       \hspace{65pt}=
        -\sum_{a=1}^d(A^a_n(t,x)-   A^a(t,x))\partial_aU^k,\\ 
        (U^{n,k}-U^k)(0,x)=0.
      \end{cases}
  \end{equation*}
  So again by Lemma \ref{lem:sec1-outline:2} and inequality  
\eqref{eq:sec2-lemmas:15}, we achieve
       \begin{equation}
         \label{eq:sec1:14}
         \Vert U^{n,k}(t)-U^{k}(t) \Vert_{H^{s-1}}^2
         \leq e^{\int_0^t a_{n,s}(\tau)d\tau} \sum_{a=1}^d\int_0^t\Vert 
           (A^a_n(\tau,\cdot)-A^a(\tau,\cdot))\partial_aU^k(\tau) 
         \Vert_{H^{s-1}}^2.
       \end{equation}
        The multiplication property \eqref{eq:mult} implies that
       \begin{equation*}
         \|(A^a_n(\tau,\cdot)-A^a(\tau,\cdot))\partial_aU^k(\tau) 
         \|_{H^{s-1}}\lesssim 
         \|A^a_n(\tau,\cdot)-A^a(\tau,\cdot)\|_{H^{s-1}}
         \|\partial_aU^k(\tau) \|_{H^{s-1}},
       \end{equation*}
       and since $U^k\in L^\infty([0,T];H^s)$, we conclude that
       $\|\partial_aU^k(\tau) \|_{H^{s-1}}$ is bounded. 
       So, by assumption \eqref{eq:sec1-outline:54} and inequality
       \eqref{eq:sec1:14},  we obtain that for any fixed $k$ 
       \begin{equation}
         \label{eq:sec:16}
         \lim_{n\to\infty}\sup_{0\leq t \leq T}
         \Vert U^{n,k}(t)-U^{k}(t) \Vert_{H^{s-1}}=0.
       \end{equation}
       So now, using \eqref{eq:sec:14}, \eqref{eq:sec:15} and
       \eqref{eq:sec:16}, we can accomplish the proof by a 
       three-$\epsilon$ argument. 
     \end{proof}
\subsection{The continuity of the flow map in the Sobolev spaces $H^s$.}
\label{sec:continuity-flow-map}
In this section, we turn  to the non-linear case and consider an
initial value problem for a quasi-linear  first order symmetric
hyperbolic systems 
\begin{equation*}
    \label{eq:Q}
    \tag{QLS}
      \begin{cases}
        \partial_t U + \sum_{a=1}^d A^a (U) \partial_a U =G(U)\\
         U(0,x)=u_{0}(x)
      \end{cases}
  \end{equation*}
  under the following assumptions:
  \begin{align}
    \label{eq:assump:1}
    & A^a(0) =0, \quad 
      A^a\in C^\infty(\setR^N;\mathcal{S}^N);\\ 
    \label{eq:assump:2}
    &G  \in C^1(\setR^N;\setR^N),\\ 
    \label{eq:assump:3}
    &   
      \begin{cases}
        \|D_x[G(u(x))-G(v(x))]\|_{H^{s-1}}\leq L 
        \|u-v\|_{H^{s-1}}, \\ 
        \text{for all} \ u,v \  \text{that belong to a bounded set }  \Omega \subset 
        H^s,
      \end{cases}
  \end{align}
  where $\mathcal{S}^N$ denotes the space of $N\times N$ symmetric matrices, and where
  $D_x$ is the derivative with respect to the spatial variable $x$
  and $L$ is a positive constant.

  The existence and uniqueness of solutions to the system \eqref{eq:Q} are
  known; see, for example, the corresponding theorems in~\cite{KATO}. 
  More precisely, if $s>\frac{d}{2}+1$, conditions
  \eqref{eq:assump:1}, \eqref{eq:assump:2}, and \eqref{eq:assump:3}
  hold, and $u_0\in H^s$, then there is a positive $T$ and a unique
  solution $U$ to \eqref{eq:Q} such that
  $U\in C^0([0,T];H^s)\cap C^1([0,T];H^{s-1})$.
  We prove that, under the above conditions, solutions to the initial
  value problem \eqref{eq:Q} depend continuously on the initial data
  in the $H^s$ norm.
  Hence, the initial value \eqref{eq:Q} is well posed in the Hadamard
  sense. 
\begin{thm}[The continuity of the flow map for a quasilinear symmetric
  hyperbolic system]
  \label{thm:1}
  Let $s>\frac{d}{2}+1$ and assume that
  \eqref{eq:assump:1}--\eqref{eq:assump:3} hold. 
  Let $u_0\in H^s$ and let
  $U(t)\in C^0([0,T];H^s)\cap C^1([0,T];H^{s-1})$ be the corresponding
  solution to \eqref{eq:Q} with initial data $u_0$. 
  If $\|u_0^n -u_0\|_{H^{s}}\to 0$, then for large $n$ the solutions
  $U^n(t)$ to \eqref{eq:Q} with initial data $u_0^n$ exist for
  $t\in[0,T]$, and moreover
  \begin{equation}
    \label{eq:sec1-outline:19}
    \lim_{n\to\infty}\sup_{0\leq t\leq T} \Vert U^n(t)-U(t) 
    \Vert_{H^{s}}=0.
  \end{equation}
\end{thm}

\begin{rem}
\label{rem:sec1-cont-flow-hs:3}
  This theorem was proved by  Kato
  \cite{KATO} using the theory of nonlinear semigroups. 
  Here we present a different approach to the proof of this theorem. 
  A proof of the continuity of the flow map of \eqref{eq:Q} appears
  also in \cite[Theorem 4.24]{Bahouri_2011}, but with zero right-hand
  side. 
  We adopted their idea to split the equations for $D_xU^n$ into two
  systems.
  As we mentioned in Subsection~\ref{sec:cont-depend-solut}, their
  proof suffers from the absence of the $H^{s-1}$ energy estimate
  which is a crucial tool for the proof.
\end{rem}

\begin{rem}
\label{rem:sec1-cont-flow-hs:4}
Theorem \ref{thm:1} remains valid also if the coefficient of the time
derivative is a positive definite matrix, the matrices~$A^a$, and the
right-hand side $G$ depends also on time and space variables.
In addition, it suffices to require that the matrices~$A^a$ belong to
$C^k$, where $k$ is the smallest integer greater than $s+1$;
however, in numerous applications those matrices are just smooth
functions of the dependent variable $U$.
For further details see Subsection \ref{sec:cosmological-context}.
\end{rem}
\begin{rem}
  \label{rem:Katos-example}
  One may ask which type of continuity of the flow map holds. 
  The answer was given by Kato
  \cite{KATO}, who  showed that for the simplest nonlinear symmetric
  hyperbolic system, namely, Burgers' equation
  \begin{equation*}
\begin{cases}  \partial_t U+ U\partial_x U= 0, \\
  U(0,x)=u_0(x),
      \end{cases}
  \end{equation*}
  the flow map cannot be H\"older continuous of order $\alpha$ for any
  $0<\alpha< 1$.
\end{rem}

\textbf{Proof outline.} 
We take a sequence $u_n$ such that $u_n\to u_0$ in $H^s$, and we show that 
$U^n$, the solution to \eqref{eq:Q} with initial data $u_n$, converges to~$U$, 
the solution to \eqref{eq:Q} with initial data $u_0$, in the $H^s$ norm.
We first show that $U^n$ converges to~$U$ in the $L^2$
norm. 
Since $U^n$ are bounded in the $H^s$ norm, the interpolation
\eqref{eq:tools:5} property implies that $U^n$ converges to $U$ in
$H^{s-1}$. The essential difficulty is to show the convergence in the $H^s$ 
norm.

In order to prove it, we formally differentiate system \eqref{eq:Q}, and
obtain a new symmetric hyperbolic system. 
The idea is to show that $DU^n\to DU$ in $H^{s-1}$ which will imply the
desired convergence. 
The linearization of the new system results in a system, in which the
coefficients $A^a\in H^s$, while the solutions belong to $H^{s-1}$. 
In this situation, we can apply the low regularity estimate.

We write $DU^n=W^n+Z^n$ and split the system into  two systems of
equations, one for each unknown.
The one for $W^n$ has a right-hand side and initial data independent of
$n$. 
So we can apply Lemma \ref{lem:sec1-outline:4} and conclude that $W^n$
converges to~$DU$ in the $H^{s-1}$ norm. 
The system for $Z^n$ has right-hand side terms and initial data which
tend to zero in the $H^{s-1}$ norm. 
In that situation, we can apply the low regularity energy estimate,
Lemma \ref{lem:sec1-outline:2}, and conclude that
$\|Z^n\|_{H^{s-1}}\to 0$.  
Having obtained the convergence of each of the systems, we achieve
that $\|U^n-U\|_{H^s}\to 0$.
\begin{proof}[Proof of Theorem \ref{thm:1}] We start with some
  observations.
  Let $U(t)$ be the solution of \eqref{eq:Q} with initial data
  $\|u_0\|_{H^s}$ and let $U^{n}(t)$ be the solution of \eqref{eq:Q} with
  initial data $\|u_0^{n}\|_{H^s}$.
  The time interval $[0,T]$ depends on $\|u_0\|_{H^s}$, while the time
  interval $[0,T_n]$ depends on $\|u_0^n\|_{H^s}$.
  Hence, since $\|u_0^n-u_0\|_{H^s}\to 0$, then $T_n\to T$, and for sufficiently large $n$ both time intervals coincide.
  Moreover, for $t\in [0,T]$ the norms
  $\|U(t)\|_{H^{s}}, \|U^n(t)\|_{H^{s}}$ are bounded by a constant $C$
  independent of $n$ and $\{U(t),U^n(t)\}$ belongs to a compact subset
  $K$ of $\setR^N$.
 
\textbf{Step 1.}
Set $V^n=U^n-U$; then it satisfies the equation
\begin{equation}
  \label{eq:proof:1}
    \begin{cases}
      \partial_t V^n + \sum_{a=1}^d A^a (U^n) \partial_a V^n
      \\
      \hspace{30pt}=G(U^n)-G(U)-\sum_{a=1}^d
      (A^a (U^n)-A^a(U))\partial_a      U,\\
       V^n(0,x)=u^n_{0}(x)-u_0(x).
    \end{cases}
\end{equation}
Therefore, by the standard $L^2$-energy estimates, we obtain
\begin{equation}
  \begin{split}
    \label{eq:sec1-cont-flow-hs:7}
    \frac{d}{dt}\|V^n(t)\|_{L^2}^2  \leq\ & a_{\infty,n}(t)
    \|V^n(t)\|_{L^2}^2 +
    \|G(U^n(t))-G(U(t))\|_{L^2}^2 \\ 
    & +\sum_{a=1}^d \|(A^a (U^n(t))-A^a(U(t)))\partial_a U(t)\|_{L^2}^2,
  \end{split}
\end{equation}
where
$a_\infty(t):= \sum_{a=1}^d\|(\partial_a  \left(A^a (U^n(t)\right)\|_{L^\infty}+1.$ 
By the Sobolev embedding theorem and the nonlinear estimate
Proposition \ref{prop:H^s}, equation (\ref{nonlinear}), it follows that
\begin{equation}
 \label{eq:proof:12}
 \|\partial_a \left(A^a U^n(t)\right)\|_{L^\infty}\lesssim 
\|\partial_a \left( A^a U^n(t)\right)\|_{H^{s-1}}\lesssim \| 
A^a(U^n(t))\|_{H^s}\lesssim \|U^n(t)\|_{H^s}.
\end{equation}
Thus $a_{\infty,n}(t)\leq C$, a constant independent of $n$.
We now estimate the difference $\| G(U^n)- G(U)\|_{L^2}^2$
following the ideas of the proof of inequality \eqref{eq:tools:6} in Proposition 
\ref{prop:H^s}. 
So we start with the expression
\begin{equation*}
 G(U^n)-G(U)=\int_0^1D_UG(\tau U^n+(1-\tau)U)(U^n-U)d\tau,
\end{equation*}
and recall that the terms $U^n(t)$, $U(t)$ are contained in the compact set $K\subset 
\setR^N$, which implies that 
\begin{equation}
  \label{eq:proof:10}
  \|G(U^n)-G(U)\|_{L^2}\leq \|D_UG\|_{L^\infty(K)}\|U^n-U\|_{L^2}.
\end{equation}
Using a similar argument, we conclude that
\begin{equation}
  \label{eq:proof:11}
  \|(A^a(U^n)-A^a(U))\partial_a U\|_{L^2}\leq 
  \| \partial_a U\|_{L^\infty}\|D_U A^a\|_{L^\infty(K)}\|U^n-U\|_{L^2}.
\end{equation}
In addition, 
$\|\partial_a U\|_{L^\infty}\lesssim \|\partial_a   
U\|_{H^{s-1}}\lesssim \| U\|_{H^{s}} $, which is bounded.
Thus inserting inequalities \eqref{eq:proof:11} and \eqref{eq:proof:10} into 
\eqref{eq:sec1-cont-flow-hs:7}  we obtain that
\begin{equation}
\label{eq:proof:14}
 \frac{d}{dt} \|V^n(t)\|_{L^2}^2 \leq C 
\|V^n(t)\|_{L^2}^2, \quad t\in [0,T],
\end{equation} 
where the constant $C$  does not depend on $n$. Applying  Gronwall's inequality 
to \eqref{eq:proof:14}, we finally arrive at the following inequality: 
\begin{equation*}
  \label{eq:proof:2}
\|V^n(t)\|_{L^2}^2\leq
    e^{Ct}\|u_0^n-u_0\|_{L^2}^2
    \leq e^{Ct}\|u_0^n-u_0\|_{H^s}^2,
\end{equation*}
and since $ \|u_0^n-u_0\|_{H^s} 
\to 0$ by the assumptions, we conclude that
\begin{equation}
  \label{eq:proof:3}
\lim_{n\to\infty}\sup_{0\leq t\leq T}\|V^n(t)\|_{L^2}^2=0.
\end{equation}
Combining it with the interpolation theorem, inequality
\eqref{eq:tools:5}, and the boundedness of $U^n$ in $H^s$, we obtain
\begin{equation}
\label{eq:proof:5}
 \sup_{0\leq t\leq T}\|V^n(t)\|_{H^{s'}}
 \lesssim \sup_{0\leq t\leq T}\|V^n(t)\|_{L^2}^{1-\frac{s'}{s}}
\end{equation} 
for any $0<s'<s$.
In particular, \eqref{eq:proof:3} and  \eqref{eq:proof:5} imply  that 
\begin{equation}
  \label{eq:sec2-main-result:6}
  \lim_{n\to\infty}\sup_{0\leq t \leq 
    T}\|U^n(t)-U(t)\|_{H^{s-1}} = 0.
\end{equation}

\textbf{Step 2.} Let
$D_xU=D_xU(t,x)=(\frac{\partial}{\partial x_j}U^i(t,x))$,
$i=1,\ldots,N$ and $j=1,\ldots, d$, be the matrix of first order derivatives.
Then by the chain rule it satisfies the system
\begin{equation}
  \label{eq:sec4-proof:5}
    \begin{cases}
       \partial_t (D_xU^n) + \sum_{a=1}^d A^a (U^n) \partial_a (D_xU^n) \\
\hspace{30pt}=      -\sum_{a=1}^d (D_UA^a)(U^n)
      D_xU^n\partial_a U^n +D_x(G(U^n)),\\
       D_x U^n(0,x)=D_xu_{0}^n(x).
    \end{cases}        
\end{equation}
We write
\begin{equation}
  \label{eq:sec4-proof:8}
    H^n =-\sum_{a=1}^d D_U (A^a)(U^n) 
    D_xU^n\partial_a U^n  \;\; \text{and} \;\; H =-\sum_{a=1}^d D_U (A^a)(U)  D_xU\partial_a U,
\end{equation}
and split the system \eqref{eq:sec4-proof:5} into two systems as follows: set 
$D_xU^n=W^n+Z^n$, where $W^n$ satisfies
\begin{equation}
\label{eq:sec4-proof:2}
    \begin{cases}
      \partial_t W^n + \sum_{a=1}^d A^a (U^n) \partial_a W^{n} = H  
+D_xG(U),\\
      W^n(0,x)=D_xu_{0}(x),
    \end{cases}        
\end{equation}
and
\begin{equation}
  \label{eq:sec4-proof:3}
    \begin{cases}
      & \partial_t Z^n + \sum_{a=1}^d A^a (U^n) \partial_a Z^{n} =
      H^n-H +
      D_x(G(U^n))- D_x(G(U)),\\
      & Z^n(0,x)=D_xu_{0}^n(x)-D_xu_0(x).
    \end{cases}        
\end{equation}
Note that $W^n$ satisfies a system with the right-hand side and initial
data independent of $n$,
while for the system of $Z^n$ these data tend to zero in the
$H^{s-1}$ norm.  

\textbf{Step 3.} We show $W^n\to D_xU$ in the $H^{s-1}$ norm
by applying Lemma \ref{lem:sec1-outline:4} to the system
\eqref{eq:sec4-proof:2}. 
So define $A^a_{n}(t,\cdot):=A^a(U^n(t,\cdot))$ and $A^a(t,\cdot):=A^a(U(t,\cdot))$. 
Recall that $A^a(\cdot)$ are smooth functions, hence by the nonlinear
estimate \eqref{nonlinear}
\begin{equation*}
  \Vert A^a_{n}(t,\cdot) \Vert_{H^{s}}=\Vert 
    A^a(U^n(t,\cdot))\Vert_{H^{s}}\lesssim
  \Vert U^n(t,\cdot) 
  \Vert_{H^{s}}\leq C
\end{equation*}
uniformly for $t\in [0,T]$, which implies that the condition
\eqref{eq:sec1-outline:14} is satisfied, and by the 
difference estimate \eqref{eq:tools:6}
\begin{equation}
  \label{eq:sec1-cont-flow-hs:10}
  \begin{aligned}
  \|A^a_n(t,\cdot)-A^a(t,\cdot)\|_{H^{s-1}}&
  =\|A^a(U^n(t,\cdot))-A^a(U(t,\cdot))\|_{H^{s-1}}\\
 & \lesssim \Vert U^n(t,\cdot)-U(t,\cdot)
  \Vert_{H^{s-1}}.
  \end{aligned}
\end{equation}
Hence \eqref{eq:sec2-main-result:6} implies that 
\begin{equation}
\label{eq:sec1-cont-flow-hs:8}
  \lim_{n\to \infty}\sup_{0\leq t\leq     T}\|A^a_n(t,\cdot)-A^a(t,\cdot)\|_{H^{s-1}}=0
\end{equation}
and thus condition \eqref{eq:sec1-outline:54} is also fulfilled.

Note that  $D_xU$ satisfies system 
\begin{equation}
  \label{eq:sec2-main-result:4}
    \begin{cases}
  \partial_t D_xU + \sum_{a=1}^d A^a (U) \partial_a D_xU = H  +D_xG(U),\\
 D_xU(0,x)=D_xu_0(x),
  \end{cases}
    \end{equation}
    and that $D_xG(U)\in H^{s-1}$. 
    This follows by setting $v\equiv 0$ in assumptions \eqref{eq:assump:3}. 
    Thus we can apply Lemma \ref{lem:sec1-outline:4} and conclude that
 \begin{equation}
   \label{eq:proof:9}
   \lim_{n\to\infty}\sup_{0 \leq t\leq T}
   \|W^n-D_xU\|_{H^{s-1}}=0.
 \end{equation}
 
\textbf{Step 4.}
 It remains to show that $Z^n\to 0$ in the $H^{s-1}$ norm. 
 By the $H^{s-1}$ energy estimate, Lemma  \ref{lem:sec1-outline:2}, we have 
 \begin{equation}
 \label{eq:proof:13}
 \begin{split}
 \Vert Z^n(t) \Vert_{H^{s-1}}^2 &\leq
   e^{\int_0^t a_{n,s}(\tau)d\tau}
   \bigg[ \Vert D_xu^n_0-D_xu_0
      \Vert_{H^{s-1}}^2\\
& \hspace{8pt}   + \int_0^t (  \Vert
          H^n(\tau)-H(\tau) \Vert_{H^{s-1}}^2 + 
        \Vert D_x[G(U^n(\tau))-G(U(\tau))] 
\Vert_{H^{s-1}}^2 
) d\tau \bigg],
  \end{split}
\end{equation}
where
$a_{n,s}(\tau):= C \sum_{a=1}^d \|A^a(U^n(\tau))\|_{H^s}+1$.
In a similar manner as in \eqref{eq:proof:12}, we conclude that
$a_{n,s}(\tau)\leq C$, a constant independent of $n$. 
We have to estimate the last two  
terms of  \eqref{eq:proof:13} by $\|U^n-U\|_{H^{s-1}}$ and $\|Z^n\|_{H^{s-1}}$.
Here we use  assumption \eqref{eq:assump:3} which allows us to
conclude that
 \begin{equation}
   \label{eq:proof:8} 
   \|D_x[G(U^n)-G(U)]\|_{H^{s-1}}\leq L
   \|U^n-U\|_{H^{s-1}}.
 \end{equation}
 We turn now to the estimate of
 $ \Vert H^n-H \Vert_{H^{s-1}} $, which involves some elaborated 
 computations whose result we present  in the following proposition,
 the proof of which we postpone for the moment.
 
\begin{prop}[Estimate of  $\Vert H^n-H \Vert_{H^{s-1}}  $]
  \label{lem:sec4-proof:1}
  Let $H^n$ and $H$ be defined by~\eqref{eq:sec4-proof:8}. Then the
  following estimate holds:
  \begin{equation}
    \label{eq:sec4-proof:16}
    \Vert H^n-H \Vert_{H^{s-1}}
    \leq C \{ \Vert U^n-U \Vert_{H^{s-1}}
+      \Vert W^n-D_xU \Vert_{H^{s-1}}
+      \Vert Z^n \Vert_{H^{s-1}}
    \}.
  \end{equation}
\end{prop}

We recall that $a_{n,s}(\tau)\leq C$ and use \eqref{eq:sec4-proof:16}, then
we insert this expression and inequality \eqref{eq:proof:8} into
inequality \eqref{eq:proof:13}, and obtain 
\begin{equation}
 \label{eq:sec4-proof:27}
\begin{split}
 \Vert Z^n(t) \Vert_{H^{s-1}}^2  \leq\ &
   e^{Ct} \bigg[ \Vert D_xu^n_0\!-\!D_xu_0
      \Vert_{H^{s-1}}^2  \\
      & +\!   C\int_0^t (  \Vert
          U^n(\tau)\!-\!U(\tau) \Vert_{H^{s-1}}^2\! +\! 
        \Vert Z^n(\tau) \Vert_{H^{s-1}}^2 \!+\!  \Vert
          W^n(\tau)\!-\!D_xU(\tau) \Vert_{H^{s-1}}^2)d\tau \bigg].
  \end{split}
\end{equation}
Thus there exists  a $T^{*}$ such that $0<T^{\ast}\leq T $  and 
\begin{equation*}
  \label{eq:sec4-proof:27a}
  \begin{split}
    \sup_{0\leq t\leq T^{\ast}} \Vert Z^n(t) &
\Vert_{H^{s-1}}^2\\
 \leq\ & 2e^{CT^{\ast}}
   \Vert D_xu^n_0-D_xu_0 \Vert_{H^{s-1}}^2 
    \\ & 
      + \Big( \sup_{0\leq t\leq T^{\ast}} \Vert U^n(t)-U(t) 
\Vert_{H^{s-1}}^2 +
        \sup_{0\leq t\leq T^{\ast}} \Vert W^n(t)-D_xU(t) 
\Vert_{H^{s-1}}^2\Big).
  \end{split}
\end{equation*}
 Now  $\Vert D_xu^n_0-D_xu_0  \Vert_{H^{s-1}}\to 0$ by the assumptions, and  with the 
combination of inequality \eqref{eq:sec2-main-result:6} and limit
\eqref{eq:proof:9} we conclude that 
\begin{equation*}
 \lim_{n\to\infty}\sup_{0\leq t\leq T^{\ast}} 
\|Z^n(t)\|_{H^{s-1}}=0.
\end{equation*}
Thus we have obtained that 
 \begin{math}
  \lim_{n\to\infty} \sup_{0\leq t\leq 
T^{\ast}}\|D_xU^n(t)-D_xU(t)\|_{H^{s-1}}=0
 \end{math},  
 and by \eqref{eq:proof:3},
  \begin{math}
  \lim_{n\to\infty} \sup_{0\leq t\leq 
T}\|U^n(t)-U(t)\|_{L^2}=0
\end{math}. Hence 
\begin{equation*}
 \lim_{n\to\infty}\sup_{0\leq t\leq T^{\ast}} 
\|U^n(t)-U(t)\|_{H^{s}}=0.
\end{equation*}

\textbf{Step 5.}
Since $T^{\ast}$ depends just on the initial data and inequalities  
\eqref{eq:proof:10} and~\eqref{eq:proof:11},   we can repeat the same 
arguments with initial data $U(T^{\ast})$ and\linebreak after the final 
iterations steps we will derive \eqref{eq:sec1-outline:19}.
This completes the proof of Theorem~\ref{thm:1}.
\end{proof}

It remains to prove Proposition \ref{lem:sec4-proof:1}.
\begin{proof}[Proof of Proposition \ref{lem:sec4-proof:1}]
  We start with
  \begin{equation*}
    \label{eq:sec4-proof:17}
    \begin{split}
       H^n-H =\ & (D_U A^a)(U^n)D_xU^n 
\partial_aU^{n}-(D_U A^a)(U)D_xU
      \partial_aU \\ 
      =\ & D_xU^n \partial_a U^n((D_U 
A^a)(U^n)-(D_U A^a)(U)) \\ 
& + (D_U A^a)(U) ( 
D_xU^n \partial_a
        U^n - D_xU \partial_aU ) =:  H_1+H_{2}.
    \end{split}
  \end{equation*}
  Now, since $A^a(\cdot)$ are smooth functions and $U$ and $U^n$ are
  bounded in the $H^{s}$ norm by a constant independent of $n$, we
  obtain by the difference estimate \eqref{eq:tools:6} and the
  multiplication property \eqref{eq:mult} that
  \begin{equation}
    \label{eq:sec4-proof:19}
    \Vert H_1  \Vert_{H^{s-1}} \lesssim \Vert D_xU^n 
    \Vert_{H^{s-1}} \Vert \partial_a U^n
    \Vert_{H^{s-1}} \Vert U^n-U \Vert_{H^{s-1}} \lesssim
     \Vert U^n-U \Vert_{H^{s-1}}.
  \end{equation}
  Next, note that
  \begin{equation*}
    \label{eq:sec4-proof:20}
    \begin{split}
      D_xU^n \partial_a U^n -D_xU \partial_a U &= ( W^n +Z^n 
)\partial_aU^n-D_xU\partial_a U \\ &= 
      Z^n \partial_aU + (W^n -D_xU)\partial_a U^n + D_xU 
(\partial_aU^n -\partial_a U).
    \end{split}
  \end{equation*}
  So by the multiplication property in $H^{s-1}$, \eqref{eq:mult:weight} and the nonlinear 
  estimate \eqref{nonlinear} in Proposition 8, we have that
\begin{equation}
  \label{eq:sec4-proof:30}
\begin{split}
     \|H_2\|_{H^{s-1}}      \lesssim\ & \|(D_U
      A^a)(U)\|_{H^{s-1}}\times \{\|Z^n\|_{H^{s-1}}
      \|\partial_a U\|_{H^{s-1}}   \\ 
 & \hspace{8pt}  + 
      \|W^n-D_xU\|_{H^{s-1}} \|\partial_a U^n\|_{H^{s-1}}
      +\|D_xU\|_{H^{s-1}} \|\partial_a U^n-\partial_a
        U\|_{H^{s-1}} \} \\ 
         \lesssim\ &
    \{\|Z^n\|_{H^{s-1}} +
      \|W^n-D_xU\|_{H^{s-1}} +\|\partial_a U^n-\partial_a
        U\|_{H^{s-1}}\}.
  \end{split}
\end{equation}
  We now express $(\partial_a U^n-\partial_a U)$ in terms of $W^n$ and 
$Z^n$, which leads to
\begin{equation*}
  \label{eq:sec2-main-result:2}
  D_x U^n = \Big( \frac{\partial}{\partial x_j}U^i  \Big)^{n}=Z^n_{ij}+W^n_{ij},
\end{equation*}
and moreover, we have that 
  \begin{equation*}
    \label{eq:sec2-main-result:5}
    ( \partial_a U^n )_{i} = Z^n_{ia}+W^n_{ia},
  \end{equation*}
  which are the rows of the split  matrix of derivatives.  Hence
  \begin{equation}
    \label{eq:proof:7}
    \begin{split}
    \|\partial_a U^n-\partial_a U\|_{H^{s-1}}
    & \leq \|W^n_{ia}-\partial_a U \|_{H^{s-1}}+\| 
Z^n_{ia}\|_{H^{s-1}} \\
    & \leq \|W^n-D_xU\|_{H^{s-1}}+\| Z^n\|_{H^{s-1}}.
    \end{split}
  \end{equation}
  
  From \eqref{eq:sec4-proof:19}, \eqref{eq:sec4-proof:30} and
  \eqref{eq:proof:7} we obtain \eqref{eq:sec4-proof:16} and that
  completes the proof of Proposition \ref{lem:sec4-proof:1}.
\end{proof}

\subsection{Sufficient conditions for $G(U)$.}
\label{subsec1:sufficient}
The Lipschitz continuity assumption\linebreak \eqref{eq:assump:3} is suitable
in   cases in which the function $G$ is a composition of a linear
and a nonlinear function.
Such a composition occurs, for example, in the case of the
Euler--Poisson system, which is discussed in detail in Section
\ref{sec:euler-poisson-makino}. 
It is interesting to present an alternative sufficient condition of
the function $G$ which is easier to check, and such that the
conclusions of Theorem \ref{thm:1} remain valid. 
We do it in the following proposition.

\begin{prop}[Sufficient condition for the source term]
  \label{prop:4}
  Theorem \ref{thm:1} remains valid if we replace the Lipschitz
  continuity of assumptions \eqref{eq:assump:3} by
  \mbox{$G\!\!\in\! C^{k}(\setR^N\!;\setR^N)$}, where $k$ is the least integer larger than
  $s+1$. 
\end{prop}

\begin{proof}
  We compute the derivative $D_x$ of the difference $G(u)-G(v)$ by the
  chain rule and obtain
  \begin{equation}
    \label{eq:suff:1}
    \begin{split}
       D_x[G(u)-G(v)]&=D_u G(u)D_x u - D_u G(v)D_x v \\ 
      & =D_uG(u)[D_x u -D_x v]+[D_u G (u)-D_u G(v)]D_x v.
    \end{split}
  \end{equation}
  Since this proposition is applied for solutions of the system
  \eqref{eq:Q}, we may assume that
  $\| u\|_{H^s}$ and $\| v\|_{H^s}$ are bounded by a positive constant $C$ and 
that
  the functions $u(x), v(x)$ are contained in a compact set $K$ of
  $\setR^N$.
  We now start with the estimate of \eqref{eq:suff:1}:
  \begin{equation}
    \label{eq:suff:2}
    \| D_u G(u)[D_x u -D_x 
        v]\|_{H^{s-1}}\lesssim \| 
D_u G(u)\|_{H^{s-1}}\| (D_x u-D_x 
v)\|_{H^{s-1}}.
  \end{equation} 
  By the nonlinear estimate, inequality \eqref{nonlinear} of
  Proposition \ref{prop:H^s}, we find that
  \begin{equation*}
    \label{eq:suff:3}
    \| D_u G(u)\|_{H^{s-1}}\leq C 
    \| u\|_{H^{s-1}},
  \end{equation*}
  where the constant $C$ depends on
  $\|D_u G\|_{C^{k-1}(K)}$ and $\| 
u\|_{L^\infty}$.  Here we need  \mbox{$k-1\geq s-1$} in order to apply inequality 
\eqref{nonlinear}.
  As to the second term of \eqref{eq:suff:1}, we have 
  \begin{equation*}
    \|[D_u G(u)-D_u G(v)]D_x 
v\|_{H^{s-1}}\lesssim
    \|D_u G(u)-D_v G(v)\|_{H^{s-1}}
    \|D_x v\|_{H^{s-1}}.
  \end{equation*}
  Applying the difference estimate \eqref{eq:tools:6} to the function
  $\frac{\partial G}{\partial u}$, we have 
  \begin{equation}
    \label{eq:suff:4}
    \|[D_u G (u)-D_u G(v)]\|_{H^{s-1}}
    \leq C \|u-v\|_{H^{s-1}},
  \end{equation}
  where $C$ depends on
  $\|D^2_u G\|_{C^{k-2}(K)} $ and 
  $\|u\|_{L^\infty}$, $\|v\|_{L^\infty}$, $\|u\|_{H^{s-1}}$,
  $\|v\|_{H^{s-1}}$. Here we need that  $k-2\geq 
s-1$, or $k\geq s+1$, in order to apply inequality \eqref{eq:tools:6}.
  Hence, we conclude from the  inequalities
  \eqref{eq:suff:2}--\eqref{eq:suff:4} that
  \begin{equation}
    \label{eq:suff:5}
    \| D_x[G(u)-G(v)]\|_{H^{s-1}}\leq C\{\|D_x u- 
        D_x v\|_{H^{s-1}}+ \| u- v\|_{H^{s-1}}\},
  \end{equation}
  whereby the Sobolev embedding theorem with constant $C$ depends just
 on the $H^s$ norm of $u $ and $v$. 
 
  So now we apply inequality \eqref{eq:suff:5} to the solutions $U^n$
  and $U$ of the system \eqref{eq:Q}, and combine it with  the splitting
  $$D_x U^n=Z^n+W^n,$$
   and obtain
  \begin{equation}
    \label{eq:suff:6}
    \begin{split}
      \| D_x[G(U^n)-&G(U)]\|_{H^{s-1}}\\ & \leq C\{
        \|U^n -U\|_{H^{s-1}}
          +    \| Z^n\|_{H^{s-1}}+\|W^n- D_x   U\|_{H^{s-1}}\}.
    \end{split}
  \end{equation}
 
  Returning to the proof of Theorem \ref{thm:1}, we replace inequality
  \eqref{eq:proof:8} by the above inequality \eqref{eq:suff:6}, which
  results in the fact that \eqref{eq:sec4-proof:27} holds under the assumptions of
  the Proposition, and hence the conclusions of Theorem \ref{thm:1}
  hold as well.
 
\end{proof}

\begin{rem} 
\label{rem:Lipschitz}
Kato's condition for Theorem \ref{thm:1} is that $G(u)$ is Lipschitz
continuous in the $H^s$ norm~\cite{KATO}. 
A sufficient condition for this property to be satisfied is that
$G\in C^{k}(\setR^N;\setR^N)$, where $k$ is the least integer larger than
$s+1$.
This follows from the difference estimate \eqref{eq:tools:6} of
Proposition \ref{prop:H^s} in
\ref{sec:mathematical-tools} (see also \cite[Theorem IV]{KATO}). 
So it turns out that our conditions and the one imposed by Kato are
not equivalent, however when applying it to the sort of systems we
consider in which $G$ is a nonlinear function of the unknowns, both
conditions lead to the same restrictions on the nonlinearities, which
means that for practical purposes there is no actual difference in
both conditions.
\end{rem}

\biblio

\section{Continuity of the flow map in $H_{s,\delta}$}
\label{sec:continuity-flow-map-1}
The ordinary Sobolev spaces $H^s(\setR^d)$  are not suited in some
settings.
One example concerns the Euler equations if the density does not have
compact support but fall off at infinity, another one regards the
Einstein equations in asymptotically flat space-times.

In those contexts, a more appropriate class of functions is represented
by the weighted Sobolev spaces in which the weights vary with the
order of the derivatives. 
These types of spaces were introduced by Nirenberg and
Walker~\cite{nirenberg73:_null_spaces_ellip_differ_operat_r} and
independently by Cantor~\cite{cantor75:_spaces_funct_condit_r}. 
Triebel extended them to fractional order and proved basic properties
such as duality, interpolation, and density of smooth
functions~\cite{triebel76:_spaces_kudrj2}. 
These types of space have also numerous applications to elliptic PDEs,
in particular to problems that arise from geometry.

Recently the authors proved the existence and uniqueness of classical
solutions for first order symmetric hyperbolic systems \eqref{eq:Q} in
these spaces \cite{BK9}.
In this section, we shall present and prove the continuity of the flow
map in these weighted spaces, and hence establish the well-posedness
in the Hadamard sense in these spaces.

\subsection{The $H_{s,\delta}$ spaces.}
\label{sec:h_s-delta-spaces}
The weighted Sobolev spaces  of integer order can be defined as completion of 
$C_0^\infty(\setR^d)$ under the norm
\begin{equation}
\label{eq:sec2-cont-flow-weighted:3}
  \| u\|_{H_{m,\delta}}^2=\sum_{|\alpha |\leq 
    m}\int_{\setR^d}(1+|x|)^{2\delta}|\partial^\alpha u(x)|^2dx.
\end{equation}
In the case $m=0$ we denote these spaces by $L_\delta^2$, that is, 
\begin{equation*}
\|u\|_{L_\delta^2}=\|(1+|x|)^\delta u\|_{L^2}.
\end{equation*}
Triebel used a dyadic decomposition to express this norm in order to
derive various properties. 
We, however, adapted it as a different definition of these spaces.

Let $\{\psi_j\}_{j=0}^\infty$ be a  dyadic partition of unity in $\setR^d$ such 
that $ \psi_j\in C_0^\infty(\setR^d)$, \mbox{$\psi_j(x)\!\geq\!0$},
$\mbox{supp}(\psi_j)\subset \{x: 2^{j-2}\!\leq\! |x|\! \leq\! 2^{j+1}\}$,
\mbox{$\psi_j(x)\!=\!1$} on \mbox{$\{x: 2^{j-1}\!\leq\! |x|\! \leq\! 2^{j}\}$} for $j=1,2,\ldots$,
$\mbox{supp}(\psi_0)\subset\{x:|x|\leq 2\}$, $\psi_0(x)=1$ on
$\{x: |x|\leq 1\}$ and\vspace{-1pt}
\begin{equation}
  \label{eq:weighted:1}
  |\partial^\alpha  \psi_j(x)|\leq  C_\alpha  2^{-|\alpha|j},\vspace{-1pt}
\end{equation}
where the constant $C_\alpha$ does not depend on $j$.

For a function $u$ that is defined  in $\setR^d$ and  $\epsilon >0$, 
$u_\epsilon$ denotes the scaling by $\epsilon$, that is,\vspace{-3pt}
\begin{equation}
  \label{eq:sec2-cont-flow-weighted:16}
  u_{\epsilon}(x)=u(\epsilon x).\vspace{-1pt}
\end{equation}
We shall basically use the scaling with $\epsilon=2^j$, where $j$ is an 
integer.

\begin{defn}[Weighted fractional Sobolev spaces]
  \label{def:weighted:3}
  Let $s,\delta\in\setR$. The weighted Sobolev space $H_{s,\delta}$ is the set of all
  tempered distributions such that the norm\vspace{-1pt}
  \begin{equation}
    \label{eq:weighted:4}
    \|u\|_{H_{s,\delta}}^2=
    \sum_{j=0}^\infty  2^{( \delta+\frac{d}{2})2j} \| (\psi_j    u)_{2^j}\|_{H^{s}}^{2}\vspace{-1pt}
  \end{equation}
  is finite.
\end{defn}

\begin{rem}
  \label{rem:equivalence}
  Triebel proved that any other dyadic sequence
  $\{\widetilde{\psi}_j\}$ that satisfies inequality \eqref{eq:weighted:1}
  results in an equivalent norm.
  Moreover, he showed that if $s\in\setN$ then the norm
  \eqref{eq:weighted:4} is equivalent to the norm
  \eqref{eq:sec2-cont-flow-weighted:3}~\cite{triebel76:_spaces_kudrj2}.
\end{rem}

In order to derive the energy estimates, we introduce an inner product in the 
weighted space \vspace{-1pt}
\begin{equation}
  \label{eq:inner-w}
  \langle u, v\rangle_{s,\delta}=\sum_{j=0}^\infty 2^{(\delta+\frac{d}{2})2j}
  \langle (\psi_j u)_{2^j}, (\psi_j v)_{2^j}\rangle_s,\vspace{-1pt}
\end{equation}
where $\langle u, v\rangle_s$ is defined by \eqref{eq:inner-product}.
Thus $H_{s,\delta}$ is a Hilbert space.

\subsection{Low regularity energy estimates in the $H_{s,\delta}$ spaces.}
\label{sec:energy-estimates-h_s}

In this subsection, we prove the analogous energy estimate to Lemma 
\ref{lem:sec1-outline:2}, which is an essential tool for the proof of the 
continuity of the flow map. 

\begin{lem}[Low regularity energy estimate in the weighted spaces]\negthickspace
  \label{lem:energy-weighted:2}
  Let \mbox{$s\!>\!\!\frac{d}{2}+1$}, $\delta\geq -\frac{d}{2}$, and assume
  $A^a \in L^\infty([0,T]; H_{s,\delta})$,
  $F \in L^\infty([0,T]; H_{s-1,\delta+1}) $,\linebreak and $u_0\in H_{s-1,\delta+1}$.
  If $U(t)\in L^\infty([0, T ];H_{s-1,\delta+1})$ is a solution to the
  initial value problem \eqref{eq:sec1-math-prel:3}, then for
  $t\in [0,T]$ the following inequality,\vspace{-1pt}
  \begin{equation}
    \label{eq:sec1-math-prel:1}
    \Vert U(t) \Vert_{H_{s-1},\delta+1}^2\leq e^{\int_0^t 
      a_{s,\delta}(\tau)d\tau} 
    \bigg( \Vert u_0 \Vert_{H_{s-1,\delta+1}}^2+ \int_0^t     \Vert 
        F(\tau,\cdot)\Vert_{H_{s-1,\delta+1}}^2 d\tau \bigg),\vspace{-1pt}
  \end{equation}
  holds, where
  $  a_{s,\delta}(\tau):= C\sum_{a=1}^d\Vert A^a(\tau)\Vert_{H_{s,\delta}}+1$.  
\end{lem}

As in Section \ref{sec:useful-lemmas}, we first prove inequality
\eqref{eq:sec1-math-prel:1} for solutions with one more degree of
regularity of the right-hand side and the initial data, and since 
$H_{s,\delta}$ are Hilbert spaces the proof of Lemma \ref{lem:energy-weighted:2} 
is accomplished by the approximation argument. 

\begin{prop}
  \label{lem:energy-weighted}
  Let $s>\frac{d}{2}+1$, $\delta\geq -\frac{d}{2}$ and assume
  $A^a \in L^\infty([0,T]; H_{s,\delta}) $,
  \mbox{$F \in L^\infty([0,T]; H_{s,\delta+1}) $}, and $u_0\in H_{s,\delta+1}$.
  If
  $$U(t)\in C^0([0,T ];H_{s,\delta+1})\cap C^1([0,T];H_{s-1,\delta+2})$$ 
  is a solution to the initial value problem
  \eqref{eq:sec1-math-prel:3}, then for $t\in [0,T]$
  \begin{equation*}
    \label{eq:sec1-math-prel:2}
    \Vert U(t) \Vert_{H_{s-1},\delta+1}^2\leq e^{\int_0^t 
      a_{s,\delta}(\tau)d\tau} 
    \bigg( \Vert u_0 \Vert_{H_{s-1,\delta+1}}^2+ \int_0^t     \Vert 
        F(\tau,\cdot)\Vert_{H_{s-1,\delta+1}}^2 d\tau \bigg),
  \end{equation*}
  where
  $ a_{s,\delta}(\tau):= C\sum_{a=1}^d\Vert A^a(\tau)\Vert_{H_{s,\delta}}+1$.  
 
\end{prop}

\begin{proof}[Proof of Proposition \ref{lem:energy-weighted}] We follow the
  same strategy as in the proof of Proposition \ref{lem:sec1:2}.
  However, the corresponding steps, such as the commutator and
  integration by parts, are more complicated due to the fact that the
  inner product~\eqref{eq:inner-w} is an infinite sum of the
  $H^s$-inner products, which possess a scaling property~\eqref{eq:sec2-cont-flow-weighted:16}.
  That is why we shall apply these tools to each of the summands of
  \eqref{eq:inner-w}.

  Since $F$ and $u_0$ belong to $H_{s,\delta+1}$, it follows from Theorem 4.3
  in
  \cite{BK9} that $\partial_t U(t)\in C([0,T];H_{s-1,\delta+1})$, and hence we can also obtain
  an analogous identity to the one given by \eqref{eq:sec1:19} in the
  weighted space and conclude that
\begin{equation}
  \label{eq:sec2-cont-flow-weighted:6}
  \begin{aligned}
\frac{1}{2}\frac{d}{dt}\|U(t)\|_{H_{s-1,\delta+1}}^2&=\frac{1}{2}\frac{d}{dt}
    \langle U(t),U(t)\rangle_{s-1,\delta+1}= \langle U(t),\partial_t U(t)\rangle_{s-1,\delta+1}\\ 
  & = -
    \sum_{a=1}^d\langle U(t),A^a\partial_a U(t)\rangle_{s-1,\delta+1}+\langle
    U(t),F(t)\rangle_{s-1,\delta+1}.
    \end{aligned}
\end{equation}

The second line of equation \eqref{eq:sec2-cont-flow-weighted:6} can
be treated by the Cauchy--Schwarz inequality which results in
  \begin{equation}
    \label{eq:ineq:5}
    \begin{split}
      | \langle U(t),F(t)\rangle_{s-1,\delta+1}|&\leq
      \|U(t)\|_{H_{s-1,\delta+1,2}}\|F(t)\|_{H_{s-1,\delta+1,2}} \\ & \lesssim \frac{1}{2}
      (\|U(t)\|_{H_{s-1,\delta+1}}^2 + \|F(t)\|_{H_{s-1,\delta+1}}^2).
    \end{split}
  \end{equation}

  We turn now to the first term in the second line of equation
  \eqref{eq:sec2-cont-flow-weighted:6}. 
To begin, we introduce some useful quantities and, for simplicity, we
  shall write $U$ instead of $U(t)$. Let
  \begin{equation}
    \label{eq:sec2-cont-flow-weighted:13}
    \begin{split}
      E(a,j,s-1) & = \langle[(\psi_j U)_{2^j}],
      [  (\psi_j(A^a\partial_a  U))_{2^j}]\rangle_{s-1} \\
      & = \langle\Lambda^{s-1}[(\psi_jU)_{2^j}],
        \Lambda^{s-1}[(\psi_j(A^a\partial_a
              U))_{2^j}]\rangle_{L^2},
    \end{split}        
  \end{equation}
where $\Lambda^{s-1}$ is defined by \eqref{eq:notation:1}.

So with the help of \eqref{eq:sec2-cont-flow-weighted:13}, the first
term of the second line of \eqref{eq:sec2-cont-flow-weighted:6} becomes
\begin{equation}
\label{eq:sec2-cont-flow-weighted:14}
    \langle U(t),A^a
    \partial_a U(t) \rangle_{s-1,\delta+1} =
   \sum_{a=1}^d \sum_{j=0}^\infty 2^{(\delta+1+\frac{d}{2})2j}E(a,j,s-1),
  \end{equation}
and therefore it  suffices to show that
\begin{equation}
  \label{eq:ineq:4}
  \sum_{j=0}^\infty 2^{(\delta+1+\frac{d}{2})2j}| 
    E(a,j,s-1)|\lesssim  
  \|A^a\|_{H_{s,\delta}}    \| U(t)\|_{H_{s-1,\delta+1}}^2.
\end{equation}

In order to compute the $H_{s,\delta}$ norm of $A^a$, we must multiply it by a 
dyadic $\psi_j$ sequence. But on the other hand, we have to see that $ 
E(a,j,s-1) $ remains unchanged. This can be achieved by the following
manipulations. Let 
\mbox{  \begin{math}
    \Psi_k=(\sum_{j=0}^\infty \psi_j)^{-1}\psi_k
  \end{math}}; then
  \begin{math}
    \sum_{k=0}^\infty \Psi_k=1
  \end{math}.
This allows us to write 
\begin{equation*}
  \label{eq:inner:2}
  \begin{split}
     E(a,j,s-1)=\ &\langle\Lambda^{s-1}[(\psi_jU)_{2^j}],
      \Lambda^{s-1}[(\psi_j(A^a\partial_a
            U))_{2^j}]\rangle_{L^2}\\
    =\ & \bigg\langle\Lambda^{s-1}[(\psi_jU)_{2^j}],
      \Lambda^{s-1}\bigg[\bigg(\psi_j\bigg(\bigg( \sum_{k=0}^\infty
              \Psi_k\bigg)A^a\partial_a
            U\bigg)\bigg)_{2^j}\bigg]\bigg\rangle_{L^2}\\ 
            = & \sum_{k=0}^\infty
    \langle\Lambda^{s-1}[(\psi_jU)_{2^j}],
      \Lambda^{s-1}[(\psi_j( \Psi_k A^a\partial_a
            U))_{2^j}]\rangle_{L^2}.
  \end{split}
\end{equation*}
Note that $\psi_j\Psi_k\not\equiv 0$ only when $j-4\leq k\leq j+4$, therefore
  \begin{equation}
\label{eq:sec2-cont-flow-weighted:7}
    E(a,j,s-1)=\sum_{k=j-4}^{j+4} E(a,j,k,s-1),
  \end{equation}
  where
  \begin{equation*}
    \label{eq:sec2-cont-flow-weighted:8}
    E(a,j,k,s-1)=\langle\Lambda^{s-1}[(\psi_jU)_{2^j}],\Lambda^{s-1}[(\psi_j( \Psi_k A^a\partial_a
            U))_{2^j}]\rangle_{L^2}.
  \end{equation*}
  After these preparations we are almost in position to apply the
  Kato--Ponce commutator estimate \eqref{eq:kato-ponce} with the
  pseudodifferential operator $\Lambda^{s-1}\partial_a$. 
  If we used ordinary unweighted spaces without weights we would just
  proceed as in the proof of Lemma \ref{lem:sec1-outline:2}; however,
  the term $( \Psi_k A^a\partial_a U(t))_{2^j}$ causes some
  complications. 
  So in order to  ``move out'' the term  $\partial_a$ from
  $( \psi_j\Psi_k A^a\partial_aU)_{2^j}$, we proceed as follows:
  \begin{equation*}
    \label{eq:energy-w:3}
    \begin{split}
      \partial_a(\psi_j\Psi_kA^aU)_{2^j} 
      & =2^j[  (\partial_a\psi_j\Psi_kA^aU)_{2^j} +
        (\psi_j\partial_a\Psi_kA^aU)_{2^j}+ (\psi_j\Psi_k(\partial_a
          A^a)U)_{2^j}] \\ 
          &\phantom{=\ } + 2^j (\psi_j\Psi_kA^a \partial_a  U)_{2^j}.
    \end{split}
  \end{equation*}
  Thus
  \begin{equation}
    \label{eq:energy-w:13}
    \begin{aligned}
    E(a,j,k,s-1) =\ &2^{-j} \langle\Lambda^{s-1}[(\psi_jU)_{2^j} ]
      ,\Lambda^{s-1}[\partial_a(\psi_j\Psi_kA^aU)_{2^j} ]
    \rangle_{L^2} \\ 
   & - \langle\Lambda^{s-1}[(\psi_jU)_{2^j} ],
      \Lambda^{s-1} [( ( \partial_a \psi_j)\Psi_kA^aU)_{2^j}
      ] \rangle_{L^2} \\
 &   - \langle\Lambda^{s-1}[(\psi_jU)_{2^j}
      ] , \Lambda^{s-1}[
        (\psi_j\partial_a\Psi_kA^aU)_{2^j}]\rangle_{L^2}    \\ 
&    - \langle\Lambda^{s-1}[(\psi_jU)_{2^j} ],
      [\Lambda^{s-1}(\psi_j\Psi_k(\partial_a A^a)U)_{2^j} ]\rangle_{L^2}.
      \end{aligned} 
  \end{equation}

  We now  consider the first term of equation
  \eqref{eq:energy-w:13},
  and make a commutation of the operator
  $\Lambda^{s-1}\partial_a$ with $(\Psi_k A^a)_{2^j}$, which results in 
  \begin{equation} 
    \label{eq:energy-w:4}
    \begin{split}
     (\Lambda^{s-1}&\partial_a)[ (\psi_j\Psi_kA^aU)_{2^j}  ]\\  
      = \ & (\Lambda^{s-1}\partial_a)[(\psi_j\Psi_kA^aU)_{2^j}
      ]-(\Psi_k A^a)_{2^j}(\Lambda^{s-1}\partial_a)[
        (\psi_jU)_{ 2^j} ]\\ 
      &  +  (\Psi_k
     A^a)_{2^j}(\Lambda^{s-1}\partial_a)[(\psi_jU)_{
          2^j}]. 
    \end{split} 
  \end{equation}
  With respect to the first two terms on the left-hand side of
  \eqref{eq:energy-w:4}, we observe that the pseudodifferential
  operator $\Lambda^{s-1}\partial_a$ belongs to the class
  $OPS_{1,0}^s$, hence by the Kato--Ponce commutator estimate
  \eqref{eq:kato-ponce} and Sobolev embedding theorem
  ($s-1>\frac{d}{2}$) we can estimate them as follows: 
  \begin{equation}
\label{eq:sec2-cont-flow-weighted:15}
    \begin{split}
      \|
(\Lambda^{s-1}&\partial_a)[(\psi_j\Psi_kA^aU)_{2^j}
        ]-(\Psi_k A^a)_{2^j}(\Lambda^{s-1}\partial_a)[
          (\psi_jU)_{ 2^j} ] \|_{L^2} \\ &\lesssim 
      \{\|D(\Psi_k
A^a)_{2^j}\|_{L^\infty}\|(\psi_jU)_{2^j}\|_{H^{s-1}}
+\|(\Psi_k A^a)_{2^j}\|_{H^s}
        \|(\psi_jU)_{2^j}\|_{L^\infty}\}\\
     & \lesssim  \|(\psi_jU)_{2^j}\|_{H^{s-1}} 
\|(\Psi_k
        A^a)_{2^j}\|_{H^s}.
    \end{split}
  \end{equation}
  For the last term of \eqref{eq:energy-w:4}, we note that
  $\Lambda^{s-1}\partial_a=\partial_a \Lambda^{s-1}$ and then, by integration by parts and the
  symmetry of the matrices $A^a$, we conclude  that
  \begin{equation*}
    \label{eq:energy-w:5}
    \begin{split}
      - 2\langle \Lambda^{s-1}[(\psi_jU)_{2^j}], &(\Psi_k
          A^a)_{2^j}(\partial_a\Lambda^{s-1})[(\psi_j^2
            U)_{2^j}]\rangle_{L^2}\\ 
        &    = \langle
        \Lambda^{s-1}[(\psi_jU)_{2^j}], \partial_a(\Psi_k
          A^a)_{2^j}\Lambda^{s-1}[(\psi_j
            U)_{2^j}]\rangle_{L^2}.
    \end{split}
  \end{equation*}
  Hence we obtain
  \begin{equation}
    \label{eq:energy-w:9}
    \begin{split}
       2|\langle \Lambda^{s-1}[(\psi_jU)_{2^j}], (\Psi_k
          A^a)_{2^j}&(\partial_a\Lambda^{s-1})[(\psi_j
            U)_{2^j}]\rangle_{L^2}|\\ \lesssim &
      \|(\psi_jU)_{2^j}\|_{H^{s-1}} ^2\|\partial_a(\Psi_k
          A^a)_{2^j}\|_{L^\infty} \\ 
          \lesssim &  \|(\psi_jU)_{2^j}\|_{H^{s-1}}
      ^2\|\partial_a(\Psi_k A^a)_{2^j}\|_{H^{s-1}} \\
      \lesssim & \|(\psi_jU)_{2^j}\|_{H^{s-1}} ^2\|(\Psi_k
          A^a)_{2^j}\|_{H^{s}}.
    \end{split}
  \end{equation}
We combine \eqref{eq:sec2-cont-flow-weighted:15} with \eqref{eq:energy-w:9} and 
see that 
  \begin{equation}
   \label{eq:energy-w:10}
   \begin{split}
   2^{-j} |2^{-j} 
\langle\Lambda^{s-1}[(\psi_jU)_{2^j} ] 
,\Lambda^{s-1}&[\partial_a(\psi_j\Psi_kA^aU)_{2^j} ] 
\rangle_{L^2} |  \\ \lesssim  & 2^{-j}
\|(\psi_jU)_{2^j}\|_{H^{s-1}}      
^2\|(\Psi_k A^a)_{2^j}\|_{H^{s}}.
\end{split}
  \end{equation} 
  
It remains to estimate the last three terms of \eqref{eq:energy-w:13}
which  can be estimated by the Cauchy--Schwarz inequality. 
First,  we recall  a well-known  property of the $H^s$ norm. If $f$ is 
a smooth function and
$\|\partial^\alpha f\|_{L^\infty}\leq K $ for all $|\alpha|\leq N$, where
$N$ is an integer greater than $s$, then $\|fu\|_{H^s} \lesssim K \|u\|_{H^s}$.

Applying it, for example, to
$f=(\partial_a\psi_j)_{2^j} $, then by \eqref{eq:weighted:1}
we obtain
\begin{math}
\|\partial^\alpha f\|_{L^\infty}\leq K
\end{math}
 for some $ K $ and all $j$, and hence
 \begin{equation*} 
   \label{eq:energy-w:6}
   \|(\partial_a \psi_j\Psi_kA^aU)_{2^j}\|_{H^{s-1}}
   \lesssim
   \|(\Psi_kA^aU)_{2^j}\|_{H^{s-1}}.
 \end{equation*}
  Similarly,
  \begin{equation*}
    \|(\psi_j\partial_a\Psi_kA^aU)_{2^j}\|_{H^{
        s-1 } }
    \lesssim
    \|(\psi_jA^aU)_{2^j}\|_{H^{s-1}}
  \end{equation*}
  and
  \begin{equation} 
    \label{eq:energy-w:8}
    \|(\psi_j\Psi_k(\partial_a A^a)U)_{2^j}\|_{H^{ s-1 } }
    \lesssim
    \|(\psi_j(\partial_aA^a)U)_{2^j}\|_{H^{s-1}}.
  \end{equation}

Finally, taking into account the equality \eqref{eq:energy-w:13} and 
inequalities  \eqref{eq:energy-w:10}--\eqref{eq:energy-w:8},
we obtain
  \begin{equation}
    \label{eq:energy-w:11}
    \begin{aligned}
   | E(a,j,k,&s-1)|  \\\lesssim
    & \|(\psi_j
U)_{2^{j}}\|_{H^{s-1}} 
    \{\|(\Psi_k A^a U)_{2^{j}}\|_{H^{s-1}} +
      \|(\psi_j A^a U)_{2^{j}}\|_{H^{s-1}} \\ &+
      \|(\psi_j (\partial_a A^a) U)_{2^{j}}\|_{H^{s-1}} +
      2^{-j}\|(\psi_jU)_{2^{j}}\|_{H^{s-1}}
      \|({\Psi_k} A^a )_{2^{j}}\|_{H^{s}}\}.
      \end{aligned}
  \end{equation}
  So now using equations \eqref{eq:sec2-cont-flow-weighted:14} and
  \eqref{eq:sec2-cont-flow-weighted:7} leads to
  \begin{equation*}
    \label{eq:energy-w:13a}
    |\langle
    U(t),A^a\partial_a U(t)\rangle_{s-1,\delta+1}|\leq \sum_{a=1}^d
    \sum_{j=0}^\infty\sum_{k=j-4}^{j+4}
    2^{(\delta+1+\frac{d}{2})2j} | E(a,j,k,s-1)|,
  \end{equation*}
and using \eqref{eq:energy-w:11} we have
\begin{equation}
  \label{eq:sec2-cont-flow-weighted:9}
  \begin{split}
  \sum_{j=0}^\infty&\sum_{k=j-4}^{j+4}
    2^{(\delta+1+\frac{d}{2})2j} | E(a,j,k,s-1)| \\  
    \lesssim\  & \sum_{j=0}^\infty\sum_{k=j-4}^{j+4} 2^{(\delta+1+\frac{d}{2})2j}
\|(\psi_j U)_{2^{j}}\|_{H^{s-1}} \\
&\hspace{31pt} \times  
    \{\|(\Psi_k A^a U)_{2^{j}}\|_{H^{s-1}} 
     + \|(\psi_j A^a U)_{2^{j}}\|_{H^{s-1}} \\ 
     &\hspace{46pt}   +  \|(\psi_j (\partial_a A^a) U)_{2^{j}}\|_{H^{s-1}} 
 +2^{-j}\|(\psi_jU)_{2^j}\|_{H^{s-1}}\|(\Psi_k
          A^a)_{2^j}\|_{H^{s}} \}. 
  \end{split}        
\end{equation}
The right-hand side of \eqref{eq:sec2-cont-flow-weighted:9} consists
of four  different terms. 
We shall estimate only two of them, since the other terms can be dealt with
in a similar fashion and is left to the reader.
For the third term, we use the Cauchy--Schwarz inequality, the
multiplication property \eqref{eq:mult:weight} and the embedding theorem
\eqref{eq:appendix-w:1}  to obtain
  \begin{align*}
       \sum_{j=0}^\infty&\sum_{k=j-4}^{j+4} 2^{(\delta+1+\frac{d}{2})2j}
      \|(\psi_jU)_{2^j}\|_{H^{s-1}}
      \|(\psi_j(\partial_aA^a)U)_{2^j}\|_{H^{s-1}} \\
      & \lesssim  \sum_{j=0}^\infty \Big(2^{(\delta+1+\frac{d}{2})j} 
      \|(\psi_j
U)_{2^j}\|_{H^{s-1}}\Big)\left(2^{(\delta+1+\frac{d}{2})j}
      \|(\psi_j(\partial_aA^a)U)_{2^j}\|_{H^{s-1}}\right) \\
     & \lesssim  \bigg(\sum_{j=0}^\infty 2^{(\delta+1+\frac{d}{2})2j} 
      \|(\psi_jU)_{2^j}\|_{H^{s-1}}^2\bigg)^{\frac{1}{2}}
      \bigg(\sum_{j=0}^\infty 2^{(\delta+1+\frac{d}{2})2j}
      \|(\psi_j(\partial_a A^a)   U)_{2^j}\|_{H^{s-1}}^2\bigg)^{\frac{1}{2}} 
    \\
      &=  \|U\|_{H_{s-1,\delta+1}} \|(\partial_a
        A^a)U\|_{H_{s-1,\delta+1}} \lesssim \|U\|_{H_{s-1,\delta+1}}^2
      \|\partial_a A^a\|_{H_{s-1,\delta+1}} \nonumber \\
   & \lesssim 
      \|U\|_{H_{s-1,\delta+1}}^2 \|A^a\|_{H_{s,\delta}}.
  \end{align*}
  
  We now turn to the last term of \eqref{eq:sec2-cont-flow-weighted:9}. 
  By the scaling properties of the $H^s$ spaces,we obtain
  \begin{equation}
    \label{eq:energy-w:12}
    \|(\Psi_k A^a)_{2^j}\|_{H^{s}}
    =\|((\Psi_k  A^a)_{2^k})_{j-k}\|_{H^{s}}\leq C(2^{j-k})
    \|(\Psi_k A^a)_{2^k}\|_{H^{s}},
  \end{equation}
  and since $j-4\leq k\leq j+k$, the constant $C(2^{j-k})$ is bounded by
  another constant that is independent of $j$ and $k$. 
  Note also that
  \begin{math}
    (\delta+1+\frac{d}{2})2-1\leq (\delta+1+\frac{d}{2})2+\delta+\frac{d}{2}
  \end{math}
  for $\delta\geq -\frac{d}{2}-1$, hence
  \begin{equation*}
    \begin{split}
      \sum_{j=0}^\infty \sum_{k=j-4}^{j+4}
      2^{(\delta+1+\frac{d}{2})2j-j}&\|(\psi_jU)_{2^j}
      \|_{H^{s-1}}^2\|(\Psi_k A^a)_{2^j}\|_{ H^{s}}\\ 
      \leq\ & \sum_{k=0}^\infty \sum_{j=k-4}^{k+4}
      2^{(\delta+1+\frac{d}{2})2j}\|(\psi_jU)_{2^j}
      \|_{H^{s-1}}^2 2^{(\delta+\frac{d}{2})j}\|(\Psi_k A^a
        )_{2^j}\|_{ H^{s}} \\
         \leq\ & \bigg(\sum_{j=0}^\infty
        \bigg(\sum_{k=j-4}^{j+4} 2^{(\delta+1+\frac{d}{2})2j}
          \|(\psi_jU)_{2^j}\|_{H^{s-1}}^2\bigg)^2\bigg)^{\frac{1}{2}} \\ 
       &   \times 
      \bigg(\sum_{j=0}^\infty \bigg( \sum_{k=j-4}^{j+4}
          2^{(\delta+\frac{d}{2})j}\|(\Psi_k A^a )_{
              2^j}\|_{H^{s}}\bigg)^2\bigg)^{\frac{1}{2}}.
    \end{split}           
  \end{equation*}
  Now, using the simple inequality
  $(\sum_j a_j^2)^{\frac{1}{2}}\leq \sum_j |a_j|$, we have
  \begin{equation*}
    \label{eq:sec2-cont-flow-weighted:10}
       \bigg(\sum_{j=0}^\infty \bigg(\sum_{k=j-4}^{j+4} 2^{(\delta+1+\frac{d}{2})2j}
          \|(\psi_jU)_{2^j}\|_{
            H^{s-1}}^2\bigg)^2\bigg)^{\frac{1}{2}} 
             \lesssim \sum_{j=0}^\infty
      2^{(\delta+1+\frac{d}{2})2j} \|(\psi_jU)_{2^j}\|_{
        H^{s-1}}^2 \\
        = \| U\|_{H_{s-1,\delta+1}}^2.
  \end{equation*}
By \eqref{eq:energy-w:12},
  \begin{equation*}
    \begin{split}
       \bigg(\sum_{j=0}^\infty \bigg( \sum_{k=j-4}^{j+4}
       2^{(\delta+\frac{d}{2})j}\|&(\Psi_k A^a )_{
            2^j}\|_{H^{s}}\bigg)^2 \bigg)^{\frac{1}{2}}\\ \lesssim\ & 
\bigg(\sum_{j=0}^\infty \bigg(
        \sum_{k=j-4}^{j+4}
        2^{(\delta+\frac{d}{2})k}2^{(\delta+\frac{d}{2})(j-k)}\|(\Psi_k
            A^a )_{2^k}\|_{H^{s}}\bigg)^2 \bigg)^{\frac{1}{2}}\\ 
\lesssim\ & \bigg(\sum_{k=0}^\infty
      2^{(\delta+\frac{d}{2})2k}\|(\Psi_k A^a
        )_{2^k}\|_{H^{s}}^2\bigg)^{\frac{1}{2}} \\
        =\ &
\|A\|_{H_{s,\delta}}. 
    \end{split}
  \end{equation*}
The remaining terms can be estimated in a similar fashion. 
So finally  we have shown  \eqref{eq:ineq:4}.
We now combine equality
\eqref{eq:sec2-cont-flow-weighted:6} with inequalities
\eqref{eq:ineq:5} and \eqref{eq:ineq:4} to obtain
\begin{equation}
\label{eq:energy-w:15}
 \frac{d}{dt}\|U(t)\|_{H_{s-1,\delta+1}}^2\leq \bigg(C\sum_{a=1}^d 
\|A^a(t)\|_{H_{s,\delta}} 
+1\bigg)\|U(t)\|_{H_{s-1,\delta+1}}^2 
+\|F(t,\cdot)\|_{H_{s-1,\delta+1}}^2.
\end{equation} 
Finally, applying the Gronwall inequality to \eqref{eq:energy-w:15} we obtain 
\eqref{eq:sec1-math-prel:1}, and that completes the proof of the proposition. 
\end{proof}

We turn now to the proof of Lemma \ref{lem:energy-weighted:2}.
The proof follows the same lines as Lemma \ref{lem:sec1-outline:2} of
Subsection \ref{sec:useful-lemmas}.
\begin{proof}[Proof of Lemma \ref{lem:energy-weighted:2}] We first refer to
  Triebel \cite{triebel76:_spaces_kudrj2}, who proved that $C_0^\infty$ is
  dense in $H_{s,\delta}$. 
  Therefore we can find sequences
  $F^k\in L^\infty([0,T];H_{s,\delta+1})$ and
  $u_0^k\in H_{s,\delta+1}$ such that
  $\|F^k(t,\cdot)-F(t,\cdot)\|_{H_{s-1,\delta+1}}\to 0$ and
  $\|u_0^k-u_0\|_{H_{s-1,\delta+1}}\to 0$ as $k\to\infty$. 
  Let now
  $U^k\in C([0,T];H_{s,\delta+1})\cap C^1([0,T];H_{s-1,\delta+2})$ be a solution to
  the initial value problem \eqref{eq:sec1-math-prel:3} with right-hand side $F^k$ and initial data $u_0^k$. 
  The existence is assured by Theorem 4.3 in  \cite{BK9}. 
  We apply now the $L_\delta^2$ energy estimate, \cite[Lemma 4.6]{BK9}
  and the Gronwall inequality to $U^k-U$, and obtain 
\begin{equation*}
  \|U^k(t)-U(t)\|_{L_{\delta+1}}^2  \leq  e^{\int_0^t 
a_\infty(\tau)d\tau}\bigg(\|u_0^k-u_0\|_{L_{\delta+1}}^2+\int_0^t
\|F^k(\tau,\cdot)-F(\tau,\cdot)\|_{L_{\delta+1}}^2d\tau\bigg),
\end{equation*}
 where
 $ a_\infty(\tau):= C\sum_{a=1}^d(\|\partial_a 
A^a\|_{L^\infty}+ \|A^a\|_{L^\infty})+1$.
Thus it follows that $U^k(t)\to U(t)$ in the $L^2_{\delta+1}$ norm. 
On the other hand, by Proposition \ref{lem:energy-weighted}, we 
have that
\begin{equation*}
  \begin{split}
    \Vert U^k(t) \Vert_{H_{s-1,\delta+1}}^2& \leq e^{\int_0^t a_{s,\delta}(\tau)d\tau}
    \bigg( \Vert u_0^k \Vert_{H_{s-1,\delta+1}}^2+ \int_0^t\Vert
        F^k(\tau,\cdot)\Vert_{H_{s-1,\delta+1}}^2 d\tau \bigg)\\
 &   \leq e^{\int_0^t a_{s,\delta}(\tau)d\tau} \bigg( \Vert u_0
      \Vert_{H_{s-1,\delta+1}}^2+1+ \int_0^t\Vert (
          F(\tau,\cdot)\Vert_{H_{s-1,\delta+1}}^2+1) d\tau \bigg).
  \end{split}
\end{equation*}
Therefore, for each $t\in [0,T]$,  $\{U^k(t)\} $ is a bounded sequence in 
$H_{s-1,\delta+1}$, which is a Hilbert space.  Thus it weakly converges, and 
the rest of the proof is the same as in Lemma \ref{lem:sec1-outline:2}.
 
 \end{proof}
 
\subsection{Continuous dependence  of the solutions on the coefficients in $H_{s,\delta}$.}
\label{sec:cont-depend-solut-1}
\begin{lem}[Continuous dependence]
  \label{lem:weight-stability}
  Let 
  $$s>\frac{d}{2}+1, \quad\delta\geq -\frac{d}{2} ,\quad
   F\in L^{\infty}([0,T], H_{s-1,\delta+1} ),\quad
 \text{and}\quad u_0\in H_{s-1,\delta+1}.$$ 
  Let $A^a $ and $ \{A_n^a\}_{n=1}^\infty $ be symmetric matrices such that
  \begin{equation}
    \label{eq:sec2-cont-flow-weighted:1}
    \Vert A^a(t,\cdot) \Vert_{H_{s,\delta}} \leq C_0, \quad    
    \Vert A^a_n(t,\cdot)  \Vert_{H_{s,\delta}}  \leq C_0, \quad t\in[0,T],   \; 
    \forall n\in\setN
  \end{equation}
  and
  \begin{equation}
    \label{eq:sec2-cont-flow-weighted:2}
    \lim_{n\to\infty} \Vert A^a_n (t,\cdot)-A^a (t,\cdot)
    \Vert_{L^\infty([0,T]; H_{s-1,\delta})} =0.
  \end{equation}
  Let $U(t), U^n(t)\in C^0([0,T]; H_{s-1,\delta+1})$ be solutions to
  \eqref{eq:sec1-math-prel:3} with coefficients $A^a$ and~$A^a_n$
  respectively and the same data $F$ and $u_0$. 
  Then 
  \begin{equation*}
    \label{eq:sec1-outline:79}
    \lim_{n\to\infty}\sup_{0\leq t \leq T}\Vert U^n(t)- U(t) 
    \Vert_{H_{s-1,\delta+1}}=0.
  \end{equation*}
\end{lem}
\begin{proof}
  The proof follows the same lines of arguments as in the corresponding
  Lemma \ref{lem:sec1-outline:4} in the $H^s$ spaces. 
  However, some differences occur because in the weighted spaces we have
  to consider also the index $\delta$ of the weights. 
  We use the approximation sequences $\{F^k\}$ and $\{u_0^k\} $ as in the
  previous Lemma, and let $U^{n,k}$ and $U^k$ be the solutions of the
  initial value problem \eqref{eq:sec1-math-prel:3} with coefficients
  $A^a_n$, $A^a$, respectively, and with right-hand side $F^k$ and initial
  data $u_0^k$. 
  Then we write the difference $U^n(t)-U(t)$ as in \eqref{eq:sec1:16}. 
  Applying Lemma \ref{lem:energy-weighted:2} to
  $(U^n(t)-U^{n,k}(t))$ we obtain 
  \begin{equation}
    \label{eq:sec2-cont:1}
    \begin{split}
      \Vert U^n&(t)-U^{n,k}(t) \Vert_{H_{s-1,\delta+1}}^2\\
      &\leq  e^{\int_0^t
        a_{n,s,\delta}(\tau)d\tau} \bigg( \Vert u_0- u_0^k \Vert_{H_{s-1,\delta+1}}^2
       + \int_0^t \Vert
          F(\tau,\cdot)-F^k(\tau,\cdot)\Vert_{H_{s-1,\delta+1}}^2
        d\tau\bigg),
    \end{split}
  \end{equation}
  where
  $ a_{n,s,\delta}(\tau):= C\sum_{a=1}^d\Vert A^a_{n}(\tau)\Vert_{H_{s,\delta}}+1$. 
  Then by \eqref{eq:sec2-cont-flow-weighted:1},
  $a_{n,s}\leq {C}$ and this constant is independent of $n$. 
  Hence inequality \eqref{eq:sec2-cont:1} implies that
  \begin{equation*}
    \label{eq:sec2-cont:2}
    \lim_{k\to\infty}\sup_{0\leq t \leq T}
    \Vert U^n(t)-U^{n,k}(t) \Vert_{H_{s-1,\delta+1}}=0,
  \end{equation*}
  uniformly in $n$.   Similarly
  \begin{equation*}
    \label{eq:sec2-cont:3}
    \lim_{k\to\infty}\sup_{0\leq t \leq T}
    \Vert U^k(t)-U(t) \Vert_{H_{s-1,\delta+1}}=0.
  \end{equation*}
  Applying again Lemma \ref{lem:energy-weighted:2} and inequality
  \eqref{eq:sec1-math-prel:1} to $(U^{n,k}(t)-U^k(t))$ we
  obtain 
  \begin{equation}
    \label{eq:sec2-cont:4}
    \begin{split}
       \| U^{n,k}(t)-U^k(t)&\|_{H_{s-1,\delta+1}}^2 \\ \leq &
      e^{\int_0^t a_{n,s,\delta}(\tau)d\tau}\sum_{a=1}^d \int_0^t \|
        (A^a_n(\tau)-A^a(\tau))\partial_a U^{k}(\tau)\|_{H_{s-1,\delta+1}}^2d
      \tau. 
    \end{split}
  \end{equation}
  It follows from Theorem 4.3 in \cite{BK9} that
  $\partial_a U^k(t)\in H_{s-1,\delta+1}$, so the multiplication property in the
  weighted spaces, Proposition \ref{prop:H^s}, (\ref{eq:mult}), implies that
  \begin{equation}
    \label{eq:sec2-cont:5}
    \| 
      (A^a_n(\tau)-A^a(\tau))\partial_a 
      U^k(\tau)\|_{H_{s-1,\delta+1}}^2\lesssim 
    \|A^a_n(\tau)-A^a(\tau)\|_{H_{s-1,\delta}}\|\partial_a 
      U^k(\tau)\|_{H_{s-1,\delta+1}}. 
  \end{equation}
  Hence we conclude by  assumption
  \eqref{eq:sec2-cont-flow-weighted:2} and inequalities
  \eqref{eq:sec2-cont:4} and \eqref{eq:sec2-cont:5} that for any fixed
  $k$
  \begin{equation*}
    \lim_{k\to\infty}\sup_{0\leq t \leq T}
    \Vert U^{n,k}(t)-U^k(t) \Vert_{H_{s-1,\delta+1}}=0.
  \end{equation*}
  The proof is completed by a three-$\epsilon$ argument.
\end{proof}

\subsection{The continuity of the flow map for Sobolev spaces
  $H_{s,\delta}$.}
\label{sec:continuity-flow-map-3}
The existence and uniqueness of solutions to the symmetric hyperbolic system \eqref{eq:Q} in the 
$H_{s,\delta}$ spaces were achieved in \cite{BK9}.
More precisely, if $u_0\in H_{s,\delta}$, then there exists a positive $T$ and a 
unique solution $U$ to \eqref{eq:Q} such that
 \mbox{$U\in C^0([0,T]; H_{s,\delta})\cap
C^1([0,T]; H_{s-1,\delta+1})$}. Here we prove the continuity 
of the solution with respect to the initial data, and  hence 
we conclude that the Cauchy problem for symmetric hyperbolic systems  is 
well-posed in a Hadamard sense in the $H_{s,\delta}$ spaces. 
The assumptions are the same as in Subsection \ref{sec:continuity-flow-map}, 
except that \eqref{eq:assump:3} is replaced by
\begin{equation}
   \label{eq:assump:4}
  \begin{cases}
    \|D_x[G(u(x))-G(v(x))]\|_{H_{s-1,\delta+1}}\leq L 
\|u-v\|_{H_{s-1,\delta}}, \\ 
\text{for all} \ u,v  \ \text{that belong to a bounded set }\  \Omega \subset 
H_{s,\delta}. \end{cases}
  \end{equation}

  \begin{thm}[The continuity of the flow map in the weighted spaces]
    \label{thm:sec1-outline:3}
    Let $s>\frac{d}{2}+1$, $\delta\geq -\frac{d}{2}$ and assume that the
    conditions \eqref{eq:assump:1}, \eqref{eq:assump:2} and
    \eqref{eq:assump:4} hold. 
    Let $u_0\in H_{s,\delta}$ and let
    $U(t)\in C^0([0,T];H_{s,\delta})\cap C^1([0,T]; H_{s-1,\delta+1})$ be the
    corresponding solution to \eqref{eq:Q} with initial data $u_0$.
    If $\|u_0^n -u_0\|_{H_{s,\delta}}\to 0$, then for large $n$ the
    solutions $U^n(t)$ to \eqref{eq:Q} with initial data $u_0^n$ exist
    for $t\in[0,T]$, and moreover
    \begin{equation}
      \label{eq:sec2-cont-flow-weighted:4}
      \lim_{n\to\infty}\sup_{0\leq t\leq T} \Vert U^n(t)-U(t) \Vert_{H_{s,\delta}}=0.
    \end{equation}
  \end{thm}

  \begin{proof}
    The proof is based on the same ideas that we used in the
    corresponding theorem for  the un-weighted space. 
    But since the rule of the weights $\delta$ is not obvious, we shall
    highlight the relevant estimates in the weighted spaces.  

    Let $U(t)$ be a solution to the initial value problem \eqref{eq:Q}
    in the interval $[0,T]$, and let $U^n(t)$ be a solution to
    \eqref{eq:Q} with initial data $u_0^n(x)$. 
    From the existence theorem \cite[Theorem 4.3]{BK9}, it follows
    that the range of the interval $[0,T]$ depends solely on the norm
    of the initial data $\|u_0\|_{H_{s,\delta}}$. 
    And since \mbox{$\|u_0^n -u_0\|_{H_{s,\delta}}\to 0$}, for $n$ sufficiently
    large the solutions $U^n(t)$ also exist in the interval $[0,T]$. 
    Furthermore, for $t\in [0,T]$ the norms
    $\|U(t)\|_{H_{s,\delta}}, \|U(t)\|_{H_{s,\delta}}$ are bounded by a constant
    $C$ independent of $n$ and $\{U(t),U^n(t)\}$ belongs to a compact
    subset $K$ of $\setR^N$.

    We set $V^n(t)=U^n(t)-U(t)$; then it satisfies equation
    \eqref{eq:proof:1}, and by $L^2_\delta$ energy estimates
    (\cite[Lemma 4.6]{BK9})
    \begin{equation}
      \label{eq:sec2-cont:7}
\begin{split}
    \frac{d}{dt} \|V^n(t)\|_{L^2_\delta}^2     \leq\ & a_{\infty(t),n}
    \|V^n(t)\|_{L^2_\delta}^2 \\
    & +\bigg\{\| G(U^n(t))-G(U(t))\|_{L_\delta^2}^2 + \sum_{a=1}^d \|
    A^a(U^n(t))-A^a(U(t))\|_{L_\delta^2}^2\bigg\},
  \end{split}
\end{equation}
    where
    $ a_{\infty,n}(t):=
    C\sum_{a=1}^d(\|\partial_a
        A^a(U(t))\|_{L^\infty}+\|A^a(U(t))\|_{L^\infty})+1$.
    Since $U(t)$ and $U^n(t)$ are contained in a compact set of $\setR^N$, 
    $\|\partial_a A^a(U(t))\|_{L^\infty}$ and $\|A^a(U(t))\|_{L^\infty}$ 
    are bounded by a constant independent of    $n$, and consequently so is
    $a_{\infty,n}(t)$.    
By a standard method of taking a difference estimate, in a way similar to 
\eqref{eq:proof:10}, it follows that
    \begin{equation}
      \label{eq:sec2-cont:8}
      \| 
        G(U^n))-G(U)\|_{L_\delta^2}\leq 
      \|D_U G\|_{L^\infty(K)}
      \|U^n-U\|_{L_\delta^2},
    \end{equation}
    and similarly 
    \begin{equation} 
      \label{eq:sec2-cont:9}
      \| (A^a(U^n-A^a(U))\partial_a U\|_{L_\delta^2}\leq \|D_U A^a\|_{L^\infty(K)}
      \|U^n-U\|_{L_\delta^2}\|\partial_a U\|_{L^\infty}. 
    \end{equation}
    The Sobolev  embedding  in the weighted space 
\eqref{eq:tools-em-w} implies that
    $$\|\partial_a U\|_{L^\infty}\lesssim \|\partial_a 
U\|_{H_{s-1,\delta+1}},$$
 and since $U\in C^0([0,T]; H_{s,\delta})\cap
C^1([0,T]; H_{s-1,\delta+1})$ this quantity is bounded. 
    Thus inserting inequalities \eqref{eq:sec2-cont:8} and
    \eqref{eq:sec2-cont:9} in \eqref{eq:sec2-cont:7} and using the
    Gronwall inequality, we obtain 
    \begin{equation}
      \label{eq:proof-w:6}
      \sup_{0\leq t\leq T}\|V^n(t)\|_{L^2_\delta}^2\leq 
      e^{C_0T}\|u^n_0-u_0\|_{L^2_\delta}^2 
      \leq e^{CT}\|u^n_0-u_0\|_{H_{s,\delta}}^2
    \end{equation}
    for some positive constant $C_0$.
    Hence, by the weighted interpolation \eqref{eq:tools:inter-w} and
    the assumptions $\|u^n_0-u_0\|_{H_{s,\delta}}\to 0$, we obtain
    \begin{equation}
      \label{eq:proof-w:1}
      \lim_{n\to\infty}\sup_{0\leq t\leq 
        T}\|U^n(t)-U(t)\|_{H_{s-1,\delta}}=0. 
    \end{equation}
    
    As in the proof of Theorem \ref{thm:1}, we now write
    $DU^n=Z^n +W^n$, where $D=D_x$ is the derivative with respect to $x$. 
    The derivatives $W^n $ and $Z^n $ satisfy equations
    \eqref{eq:sec4-proof:2} and \eqref{eq:sec4-proof:3} respectively.
    Let $A_n^a(t,\cdot)=A^a(U^n(t,\cdot))$ and
    $A^a(t,\cdot)=A^a(U(t,\cdot))$. Then by the nonlinear estimate in the
    weighted spaces, \eqref{eq:14}, Proposition \ref{prop:w},  we obtain
    \begin{equation*}
      \|A_n^a(t,\cdot)\|_{H_{s,\delta}}=\|A^a(U^n(t,\cdot))\|_{
        H_{s,\delta }}\lesssim 
      \|U^n(t,\cdot)\|_{H_{s,\delta}},
    \end{equation*}
    and using the difference estimate in Proposition \ref{prop:w}, \eqref{eq:tools-w:2},
    we conclude that 
    \begin{equation*}
    \begin{aligned}
      \|A_n^a(t,\cdot)-A^a(t,\cdot)\|_{H_{s-1,\delta}}&=
      \|A^a(U^n(t,\cdot))-A^a(U(t,\cdot))\|_{H_{s-1,\delta}}\\
      &\lesssim 
      \|U^n(t,\cdot)-U(t,\cdot)\|_{H_{s-1,\delta}}.
      \end{aligned}
    \end{equation*}
    Therefore, \eqref{eq:proof-w:1} implies that conditions
    \eqref{eq:sec2-cont-flow-weighted:1} and
    \eqref{eq:sec2-cont-flow-weighted:2} of Lemma
    \ref{lem:weight-stability} are satisfied. 
    Also, it follows from condition \eqref{eq:assump:4} that
    $DG(U(t))\in H_{s-1,\delta+1}$. 
    Thus Lemma \ref{lem:weight-stability} implies that
    \begin{equation}
      \label{eq:proof-w:2}
      \lim_{n\to\infty}\sup_{0\leq t\leq T}\|
        W^n-DU\|_{H_{s-1,\delta+1}}=0.
    \end{equation}
    We now show that
    $\| Z^n(t)\|_{H_{s-1,\delta+1}}\to 0$.
    By assumption \eqref{eq:assump:4},
    \begin{equation}
      \label{eq:proof-w:3}
      \|D_x[G(U^n(t))-G(U(t))]\|_{H_{s-1,\delta}}
      \lesssim\|U^n(t)-U(t)\|_{H_{s-1,\delta}}.
    \end{equation}
    Replacing the $H^s$ calculus with $H_{s,\delta}$ in the proof of
    Proposition \ref{lem:sec4-proof:1}, we obtain
    \begin{equation}
      \label{eq:proof-w:4}
      \Vert H^n-H \Vert_{H_{s-1,\delta+1}}
      \lesssim \{ \Vert U^n-U \Vert_{H_{s-1,\delta}}
        +      \Vert W^n-DU \Vert_{H_{s-1,\delta+1}}
        +      \Vert Z^n \Vert_{H_{s-1,\delta+1}}   \}.
    \end{equation}
    So by the low regularity energy estimate, Lemma
    \ref{lem:energy-weighted:2}, and inequalities \eqref{eq:proof-w:3}
    and \eqref{eq:proof-w:4}, we can conclude that 
    \begin{equation*}
      \begin{split}
        \| Z^n(t)\|_{H_{s-1,\delta+1}}^2  \leq\ &e^{\int_0^t
          a_{n,s,\delta}(\tau) d\tau}\bigg[\| D_xu_0^n -
            D_xu_0\|_{H_{s-1,\delta+1}}^2 + C \int_0^t(
            \| U^n(\tau)-U(\tau)\|_{H_{s-1,\delta}}^2
        \\ 
        & +\| Z^n(\tau)\|_{H_{s-1,\delta+1}}^2+ \|
              W^n(\tau)-D_xU(\tau)\|_{H_{s-1,\delta+1}}^2)d
          \tau\bigg],
      \end{split}
    \end{equation*}
    where $ a_{n,s,\delta}:= C\sum_{a=1}^d\| A^a(U^n)\|_{H_{s,\delta}} +1$. 
    Since $a_{n,s,\delta}$ is bounded by a positive constant $C_0$ which is
    independent of $n$, there exists a $T^{*}$ such that 
    $0\leq T^\ast\leq T$ and 
    \begin{equation*}
      \begin{split}
        \sup_{0\leq t\leq T^\ast} \|          Z^n(t)\|_{H_{s-1,\delta+1}}^2 
          &\leq \ 2 e^{C_0
          T^\ast}\Big[\| D_xu_0^n -
            D_xu_0\|_{H_{s-1,\delta+1}}^2 
        \\ &\hspace{25pt}
          +CT^\ast \sup_{0\leq \tau \leq T^\ast}\{\|
              U^n(\tau)-U(\tau)\|_{H_{s-1,\delta}}^2+ \|
              W^n(\tau)-D_xU(\tau)\|_{H_{s-1,\delta+1}}^2\}\Big].
      \end{split}
    \end{equation*}
    Thus the limits \eqref{eq:proof-w:1} and \eqref{eq:proof-w:2}
    imply that
    \begin{equation*}
      \lim_{n\to\infty}\sup_{0\leq t\leq T^{\ast}}\| Z^n(t)\|_{H_{s-1,\delta+1}}=0,
    \end{equation*}
    and consequently
    \begin{equation}
      \label{eq:proof-w:5}
      \lim_{n\to\infty}\sup_{0\leq t\leq T^{\ast}}\|
        D_xU^n(t)-D_xU(t)\|_{H_{s-1,\delta+1}}=0.
    \end{equation}
    We recall that for any $u\in H_{s,\delta}$,
    $\|u\|_{H_{s,\delta}}^2\sim \|u\|_{L_\delta^2}^2+\| Du\|_{H_{s-1,\delta+1}}^2$ (see,
    e.g., \cite{triebel76:_spaces_kudrj2}), so we conclude from
    \eqref{eq:proof-w:6} and \eqref{eq:proof-w:5} that
    \begin{equation*}
      \label{eq:proof-w:7}
      \lim_{n\to\infty}\sup_{0\leq t\leq T^{\ast}}\|
        U^n(t)-U(t)\|_{H_{s,\delta}}=0.
    \end{equation*}
    The extension of the above limit to the interval $[0,T]$ is
    derived precisely as in the proof of Theorem \ref{thm:1}, which
    proves \eqref{eq:sec2-cont-flow-weighted:4}.
  \end{proof}

\begin{rem}
  The conclusions of Theorem \ref{thm:sec1-outline:3} remain valid
  even if we replace assumption \eqref{eq:assump:4} by
  $G\in C^k(\setR^N;\setR^N)$, where $k$ is the least integer greater than
  $s+1$. 
  The proof is almost identical to that of Proposition \ref{prop:4},
  except that we replace the $H^s$ calculus by  $H_{s,\delta}$.   
\end{rem}

\biblio

\section{Applications to the Euler--Poisson--Makino system}
\label{sec:euler-poisson-makino}

The Euler--Poisson system is given by
\begin{align}
  \label{eq:EP:0}
  \partial_t\rho +\sum_{a=1}^3( v^a\partial_a\rho +
  \rho \partial_a v^a)  &= 0, \\
  \label{eq:sec3-applications:5}
  \rho \bigg(\partial_t v^a +\sum_{b=1}^3 v^b\partial_b v^a \bigg) + \partial^a p &= -\partial^a\phi,\\
  \label{eq:sec3-applications:6}
  \Delta\phi &= 4\pi \rho,
\end{align}
where $\rho$ is the density and $(v^1,v^2,v^3)$ is the velocity vector field.
This system describes the motion of a gas under a self-gravitational
force. 
It consists of hyperbolic evolution equations \eqref{eq:EP:0} and
\eqref{eq:sec3-applications:5} that are coupled to the linear elliptic Poisson
equation \eqref{eq:sec3-applications:6}.
We consider this system for $x\in M$, where $M=\setR^3$ or
$M=\setT^3$, respectively, and $t\geq 0$. 
We also use the notation $\partial^a\phi=\sum_{k=1}^3\delta^{ak}\partial_k\phi$.
Moreover, we consider mostly the barotropic equation of state of the
form
\begin{equation}
  \label{eq:state}
  p=K\rho^\gamma,\quad 1<\gamma\leq3,\  0<K.
\end{equation}

We will consider the following three models to that system:
\begin{compactenum}[\upshape (1)]
  \item In the first case, the density $\rho$ has compact support. 
  Such a situation would correspond to the time evolution of a compact
  body, like a star. 
  However, the Euler equations degenerate when the density tends to
  zero. 
  This difficulty was somehow circumvented by Makino by introducing a
  new matter variable~$w$ which is a nonlinear function of the
  density $\rho$.
  Since the latter has compact support, the Sobolev spaces $H^s(\setR^{3})$  can be
  used.
  \item In the second case, the density falls off in an appropriate
  sense but could become zero at spatial infinity.
  So again a regularizing variable is used, but now the Poisson
  equation with such a source term requires a different functional
  setting; one possibility is to use weighted Sobolev spaces
  $H_{s,\delta}(\setR^3)$, which we have introduced in Subsection \ref{sec:h_s-delta-spaces}.
  \item The third and last case corresponds to  situations in
  which the density is spread all over the whole space. 
  Physically, that would correspond to a cosmological situation. 
  In that case, no degeneration of the Euler equation takes place and
  one can use the density $\rho$ as an unknown and  $H^s(\setT^3)$ in
  the functional setting. 
\end{compactenum}

The system of evolution equations \eqref{eq:EP:0} and
\eqref{eq:sec3-applications:5} is hyperbolic but not symmetric
hyperbolic without further manipulations.
It can be cast in such a form easily by choosing an appropriate
multiplier
\cite[\S 1.2]{majda84:_compr_fluid_flow_system_conser} 
leading to a system of the form 
\begin{equation*}
 A^0(U)\partial_{t}U+\sum_{k=1}^3A^k(U)\partial_kU=G(U).
\end{equation*}
In the context of the evolution of a gas which describes an isolated or a
quasi-isolated body (with appropriate fall-off condition), the density
$\rho$ is not strictly positive. 
This causes the symmetrized form of the Euler system to degenerate
since the matrix $A^0$ is then no longer positive definite. 
The only known method to circumvent this difficulty is to regularize the
system by introducing a new matter variable
\cite{makino_86}, which we will briefly discuss in  Subsections
\ref{sec:case-comp-supp} and \ref{sec:case-density-which}.

\subsection{The model of a compactly supported density.}
\label{sec:case-comp-supp}
It was observed by Makino
\cite{makino_86} that the difficulty mentioned in the previous paragraph can be, to
some extent, circumvented by using a new matter variable $w$ in place
of the density. 
For this reason, we introduce the quantity
\begin{equation*}
  \label{eq:proto-euler-poisson-banach:1}
  w=\frac{2\sqrt{K \gamma}}{\gamma-1}\rho^{\frac{\gamma-1}{2}},
\end{equation*}
which allows treating  the situation where $\rho=0$. 
Replacing the density $\rho$ by the Makino variable $w$, the system
\eqref{eq:EP:0} coupled
with the equation of state \eqref{eq:state} takes the following form:
\begin{align}
  \label{eq:iteration-scheme:1b}
  \partial_tw +\sum_{a=1}^3\Big( v^a\partial_aw + \frac{\gamma-1}{2}w 
  \partial_a v^a\Big)
  &= 0 \\
  \label{eq:iteration-scheme:2b}
  \partial_t v^a +\sum_{b=1}^3 v^b\partial_b v^a  + 
  \frac{\gamma-1}{2}w\partial^a
  w &=  -          \partial^a\phi  \\
  \label{eq:iteration-scheme:3b}
  \Delta\phi &= c_{K,\gamma}w^{\frac{2}{\gamma-1}},  
\end{align}
where
$c_{K\gamma}=
4\pi(\frac{\gamma-1}{2\sqrt{K\gamma}})^{\frac{2}{\gamma-1}}$.
We call the system
\eqref{eq:iteration-scheme:1b}--\eqref{eq:iteration-scheme:3b} the
Euler--Poisson--Makino system. 

Here we shall apply Theorem \ref{thm:1} to establish the
continuity of the Euler--Poisson--Makino system with respect to the
initial data, and hence to conclude that this system is well-posed in
the Hadamard sense.

We consider this first order symmetric hyperbolic system
\eqref{eq:iteration-scheme:1b}--\eqref{eq:iteration-scheme:2b},
coupled with the Poisson equation \eqref{eq:iteration-scheme:3b}, and
with the initial data
\begin{equation}
  \label{eq:ep:4}
    \begin{cases}
      w(0,x)=w_0(x)\geq 0,\\
      v^a(0,x)=v^a_0(x), & a=1,2,3.
    \end{cases}
\end{equation}
Existence and uniqueness of solutions to the system
\eqref{eq:iteration-scheme:1b}--\eqref{eq:iteration-scheme:3b} for
initial data \eqref{eq:ep:4} with compact support were proved by
Makino \cite{makino_86} for $1<\gamma\leq \frac{5}{3}$.
Here we apply the main result, Theorem \ref{thm:1}, and we establish
the continuity of the flow map for the system
\eqref{eq:iteration-scheme:1b}--\eqref{eq:iteration-scheme:3b}.

\begin{thm}[Continuity of the flow map of the Euler--Poisson--Makino
  system with compact density]
\label{thr:sec3-applications:2}
  Let $\frac{5}{2}<s$ when $\frac{2}{\gamma-1}$ is an integer,
  $\frac{5}{2}<s<\frac{2}{\gamma-1}+\frac{1}{2}$ otherwise. 
  Let $(w,v^a)\in C^{0}([0,T];H^s)\cap C^1([0,T];H^{s-1})$ be the solution of
  the Euler--Poisson--Makino system
  \eqref{eq:iteration-scheme:1b}--\eqref{eq:iteration-scheme:3b} with
  initial data $(w_0,v^a_0)\in H^s$, where $w_0$ has compact support. 
  If
  \begin{equation*}
    \lim_{n\to\infty}\|(w_0^n,(v_0^a)^n)-(w_0,v_0^a)\|_{H^s}\to 0,
  \end{equation*}
  then for sufficiently large $n$ the solution $(w^n, (v^a)^n)$ to
  \eqref{eq:iteration-scheme:1b}--\eqref{eq:iteration-scheme:3b} with
  initial data $( w_0^n,( v^a _{^0})^{n} )$ exists
  in the interval $[0,T]$,  and moreover,
  \begin{equation*}
    \label{eq:sec3-applications:8}
    \lim_{n\to\infty}\sup_{0\leq 
      t\leq T}\|(w^n,(v^a)^n)(t)-(w,v^a)(t)\|_{ H^s } =0.
  \end{equation*}
\end{thm}

\begin{rem}
\label{rem:sec3-applications:1}
  Makino proved local existence and uniqueness for
  $1<\gamma\leq \frac{5}{3}$, but with the combinations of the power
  estimate \eqref{eq:tools:7} and calculus in $H^s$ (Proposition
  \ref{prop:H^s}) one easily extends his result for
  $\frac{5}{2}<s<\frac{2}{\gamma-1}+\frac{3}{2}$ for a non-integer
  $s$ and for $\frac{5}{2}<s$ otherwise. 
\end{rem}

\begin{rem}
  We note that the ranges for $s$ and $\gamma$ for the existence and
  uniqueness theorems are larger than in the case of the continuity of
  the flow map. 
  Remark \ref{rem:sec3-applications:1} implies $1<\gamma\leq 3$ for
  the existence and uniqueness; however, for the continuity of the flow map 
  need we need $1<\gamma\leq 2$.
  That phenomenon is caused by the different requirements of the
  function $G$ on the right-hand side of \eqref{eq:Q}. 
  For existence, it suffices to demand that on the right-hand side of
  \eqref{eq:Q}, $G$ is  bounded in the~$H^s$ topology, 
  while for the continuity of the flow map we need $G$ to be
  Lipschitz, either in $H^s$ as Kato requires,
  \cite{KATO} or $DG$ in $H^{s-1}$ as we demand.
\end{rem}

\begin{proof}
  From the existence theorem, we know that
  $\|(w^n,(v^a)^n)(t)\|_{H^s}$ and\linebreak
  $\|(w,v^a)(t)\|_{H^s}$ are bounded by a constant $C$
  independent of $n$ for $t\in[0,T]$.
  Let $U$ denote the unknown $U=U(w,v^1,v^2,v^3)$.
  Obviously, assumption \eqref{eq:assump:1} is satisfied, so in order
  to apply Theorem \ref{thm:1} all we have to do is to check that the
  right-hand of \eqref{eq:iteration-scheme:2b} is $C^1$ and satisfies
  \eqref{eq:assump:3}, that is, $D_x[G(U)-G(\widehat U)]$ is Lipschitz
  in $H^{s-1}$ when $U$ and $\widehat U$ belong to a bounded set of
  $H^s$. 
  Now we observe that $G(U)=(0,-\nabla\phi)$, where $\phi$ is a
  solution to the Poisson equation \eqref{eq:iteration-scheme:3b}.
  Consequently (omitting the zero component)
  $G(U)= -\nabla
  \Delta^{-1}(c_{K\gamma}w^{\frac{2}{\gamma-1}})$, and
  \begin{equation*}
    \label{eq:sec3-applications:10}
    \begin{split}
      D_x[G(U)-G(\widehat U)] & =D_x [\nabla
        \Delta^{-1}(c_{K\gamma}w^{\frac{2}{\gamma-1}})-\nabla
        \Delta^{-1}(c_{K\gamma}\widehat
          w^{\frac{2}{\gamma-1}})] \\ & =
      c_{K\gamma}D_x\nabla\Delta^{-1}(w^{\frac{2}{\gamma-1}}-\widehat
        w^{\frac{2}{ \gamma-1}}),
    \end{split}
  \end{equation*}
  where $w,\widehat w\in H^s$. 
  Note that $D_x\nabla\Delta^{-1}$ is a pseudodifferential operator of
  order zero (a singular integral), which is a bounded linear operator
  on $H^s$ for all $s$ (see, e.g., 
  \cite{stein70:_singul_integ_differ_proper_funct}). 
  Therefore
  \begin{equation*}
    \begin{split}
       \| D_x[G(U)-&G(\widehat U)]\|_{H^{s-1}}
      \\ = c_{K\gamma} &
      \|D_x\nabla\Delta^{-1}(w^{\frac{2}{\gamma-1}}-\widehat
          w^{\frac{ 2}{\gamma-1}})\|_{H^{s-1}} \lesssim
      \|w^{\frac{2}{\gamma-1}}-\widehat
        w^{\frac{2}{\gamma-1}}\|_{H^{s-1}}.
    \end{split}
  \end{equation*}
 Hence it remains to estimate the power difference
  $
  \|w^{\frac{2}{\gamma-1}}- \widehat w^{\frac{2}{\gamma-1}}\|_{H^{s-1}}.$
  So setting $\beta=\frac{2}{\gamma-1}$, we have
  \begin{equation}
  \label{eq:sec3-diff}
    w^\beta-\widehat w^\beta=\beta\int_0^1(\tau 
      w^{\beta-1}+(1-\tau)\widehat w^{\beta-1})(w-\widehat w)d\tau, 
  \end{equation}
  and then by the multiplication property \eqref{eq:mult} we
  obtain
  \begin{equation*}
    \|w^\beta -\widehat w^\beta\|_{H^{s-1}}\lesssim
    \beta(\|w^{\beta-1}\|_{H^{s-1}}+ 
      \|\widehat
        w^{\beta-1}\|_{H^{s-1}})\|w-\widehat w\|_{H^{s-1}}.
  \end{equation*}
  By the fraction power estimate \eqref{eq:tools:7}, we obtain
  \begin{equation*}
    \|w^{\beta-1}\|_{H^{s-1}}\lesssim \| w\|_{H^{s-1}}\lesssim \|w\|_{H^{s}}, 
    \quad\frac{3}{2}<s-1<\beta-\frac{1}{2}, 
  \end{equation*}
  and the same estimate for $\widehat w$.
  This implies that
  \begin{equation*}
    \|w^\beta-\widehat w^\beta\|_{H^{s-1}}\lesssim
    \beta(\|w\|_{H^{s}}+ 
      \|\widehat w\|_{H^{s}})\|w- \widehat w\|_{H^{s-1}}.
  \end{equation*}
  Since $\|w\|_{H^s}, \|\widehat w\|_{H^s}\leq C$, $D_x(G(U))$ is
  Lipschitz if $\frac{3}{2}<s-1<\frac{2}{\gamma-1}+\frac{1}{2}$, which
  implies $\gamma<2$.
  For $\gamma=2$, $w^{\frac{2}{\gamma-1}}=w^{2}$, and then is obviously Lipschitz. 
\end{proof}
 
\subsection{The model  of a density which falls off at infinity.}
\label{sec:case-density-which}
The setting of the equations is the same  one as of in Section
\ref{sec:case-comp-supp}. However, since the density does not have
compact support but falls off at infinity, a different functional
setting must be used. Recently, the authors
\cite{BK9} have proven a local existence and uniqueness theorem for the 
Euler--Poisson--Makino system 
\eqref{eq:iteration-scheme:1b}--\eqref{eq:iteration-scheme:3b}, where the 
solution
\begin{equation*}
 (w,v^a)\in C^{0}([0,T];H_{s,\delta})\cap 
C^1([0,T];H_{s-1,\delta+1}), 
\end{equation*}
for $-\frac{3}{2}+\frac{2}{[\frac{2}{\gamma-1}]-1}\leq 
\delta<-\frac{1}{2}$,
$ \frac{5}{2}<s$ if $\frac{2}{\gamma-1}$ is an integer, and
$\frac{5}{2}<s<\frac{5}{2}+\frac{2}{\gamma-1}-[\frac{2}{\gamma-1}]$
otherwise. Here $[\cdot ]$ denotes the integer part of a real number.

Those conditions restrict $\gamma\in (1,\frac{5}{3})$.
We shall prove the continuity of the flow map under the same bounds of the 
parameters $\delta$, $s$, and $\gamma$. 

\begin{thm}[Continuity of the flow map of the Euler--Poisson--Makino\linebreak
  system with density which falls off at infinity]\
  \label{thr:sec3-applications:1}
  Let 
  $1<\gamma < {5}/{3}$,\linebreak
  $-\frac{3}{2}+\frac{2}{[\frac{2}{\gamma-1}]-1}\!\leq\!
  \delta\!<\!-\frac{1}{2}$, $ \frac{5}{2}\!<\!s$ if $\frac{2}{\gamma-1}$ is an integer
  and
  $\frac{5}{2}\!<\!s\!<\!\frac{5}{2}+\frac{2}{\gamma-1}-[\frac{2}{\gamma-1}]$
  otherwise. 
  Let \mbox{\begin{math}
    (w,v^a)\in C^{0}([0,T]; H_{s,\delta})\cap C^1([0,T];
      H_{s-1,\delta+1})
  \end{math}}
  be the solution of the Euler--Poisson--Makino system
  \eqref{eq:iteration-scheme:1b}--\eqref{eq:iteration-scheme:3b}
  with initial data $(w_0, v^a_0)\in H_{s,\delta}$ and $w_0\geq 0$.
  If
  \begin{equation*}    
\lim_{n\to\infty}\|(w_0^n,(v_0^a)^n)-(w_0,v_0^a)\|_{H_{s,\delta}}\to 
0,
  \end{equation*}
  then for sufficiently large $n$ the solution $(w^n, (v^a)^n)$ to
  \eqref{eq:iteration-scheme:1b}--\eqref{eq:iteration-scheme:3b} with
  initial data $(w_0^n,(v_0^a)^n)\in H_{s,\delta}$ exists in the interval
  $[0,T]$, moreover,
  \begin{equation*}
    \label{eq:sec3-applications:2}    
    \lim_{n\to\infty}\sup_{0\leq 
      t\leq T}\|(w^n,(v^a)^n)(t)-(w,v^a)(t)\|_{H_{s,\delta } }= 0.
  \end{equation*}
 
\end{thm}

\begin{proof}
  The proof has some similarities to that of Theorem
  \ref{thr:sec3-applications:2}, but there are some technical
  difficulties due to the use of weighted Sobolev spaces
  $H_{s,\delta}$.
  We have to verify that $G(U)$ satisfies condition
  \eqref{eq:assump:4}, and then we apply Theorem
  \ref{thm:sec1-outline:3}.
  As in the previous section, $G(U)=-(0,\nabla \phi)$, where
  $\Delta \phi=c_{K\gamma}w^{\frac{2}{\gamma-1}}$. 
  We shall use the embedding property \eqref{eq:appendix-w:1} and the
  commutator of the operators $\nabla\Delta^{-1}$. 
  Then, setting $\beta=\frac{2}{\gamma-1}$, we have
\begin{align*}
  \|D_x[G(U)-G(\widehat U)]\|_{H_{s-1,\delta+1}}  
  &  = c_{K\gamma}
    \|D_x[\nabla \Delta^{-1}w^{\beta}-\nabla 
    \Delta^{-1}\widehat w^{\beta}]\|_{H_{s-1,\delta+1}}\\
  & \lesssim
    \|\nabla \Delta^{-1}( w^{\beta} - \widehat w^{\beta} 
    )\|_{H_{s,\delta}}
    =  
    \| \Delta^{-1}\nabla( w^{\beta}-\widehat w^{\beta} 
    )\|_{H_{s,\delta}}.
\end{align*}
  Now we use the fact that the Laplace operator
  \begin{equation*}
    \Delta:H_{s,\delta}\to H_{s-2,\delta+2}
  \end{equation*}
  is an isomorphism for $-\frac{3}{2}<\delta<-\frac{1}{2}$;
  see \cite{BK7} or
  \cite{maxwell06:_rough_einst}.  That implies
  \begin{equation*}
    \label{eq:sec3-applications:12}
    \| \Delta^{-1}\nabla( w^{\beta}-\widehat w^{\beta}
    )\|_{H_{s,\delta}}
    \lesssim 
      \| \nabla( w^{\beta}-\widehat w^{\beta}
      )\|_{H_{s-2,\delta+2}}
    \lesssim \| ( w^{\beta}-\widehat w^{\beta} 
)\|_{H_{s-1,\delta+1}}.
  \end{equation*}
  Hence it remains to estimate the difference
  $ \|w^{\beta}-\widehat w^{\beta}\|_{H_{s-1,\delta+1}}$. 
 
  We use again the identity \eqref{eq:sec3-diff} and the multiplicity
  property \eqref{eq:mult:weight} of Proposition~\ref{prop:H^s},
  the nonlinear difference estimate, and obtain
  \begin{equation*}
    \|w^\beta-\widehat w^\beta\|_{H_{s-1,\delta+1}}\lesssim
    \beta(\|w^{\beta-1}\|_{H_{s-1\delta+1}}+ 
      \|\widehat w^{\beta-1}\|_{H_{s-1,\delta+1}})\|w-\widehat w\|_{H_{s-1,\delta}}.
  \end{equation*}
  
  Thus if we estimate $\|w^{\beta-1}\|_{H_{s-1,\delta+1}}$ by
  $\|w\|_{H_{s,\delta}}$, then the proof is completed precisely as
  in the proof of Theorem \ref{thr:sec3-applications:2}. 
  This will be done by the following proposition.
\end{proof}

\begin{prop}[Nonlinear estimate of the power of functions]
  \label{prop:power}
  Suppose that $ w\in H_{s,\delta}$, $0\leq w$ and $\beta$ is a real number greater than 
  or equal to $3$. 
  Then:
  \begin{compactenum}[\upshape (1)]
    \item If $\beta$ is an integer, $\frac{3}{2}<s$, and
    $\frac{1}{\beta-2}-\frac{3}{2}\leq \delta$, then
    \begin{equation}
      \label{eq:non:1}
      \|w^{\beta-1}\|_{H_{s-1,\delta+1}}\leq C 
      (\|w\|_{H_{s,\delta}})^{\beta-1}.
    \end{equation}
    \item If $\beta\not\in \setN$,
    $\frac{5}{2}<s<\beta-[\beta]+\frac{5}{2}$, and
    $\frac{1}{[\beta]-2}-\frac{3}{2}\leq \delta$, then
    \begin{equation*}
      \label{eq:non:2}
      \|w^{\beta-1}\|_{H_{s-1,\delta+1}}\leq C 
      (\|w\|_{H_{s,\delta}})^{[\beta]-1}.
    \end{equation*}
  \end{compactenum}
\end{prop}

\begin{proof}
  We consider the case $\beta\in \setN$ and apply equation
  \eqref{eq:tools-w:1} of Lemma \ref{lem:improved-mult}, with $u_i=w$,
  $i=1,\ldots,\beta-1$. 
  This requires that
  $(\delta+1)\leq (\beta-1)\delta+(\beta-2)\frac{3}{2}$, or equivalently
  $\frac{1}{\beta-2}-\frac{3}{2}\leq \delta$, and $s>\frac{3}{2}$. 
  Hence we get \eqref{eq:non:1}. 
  For $\beta\not\in \setN$ it can be proven in a similar fashion as in
  \cite[Prop.~4.10]{BK9}, using Kateb's power estimate \eqref{eq:kateb}.
\end{proof}

Recall that for the isomorphism of $\Delta: H_{s,\delta}\to H_{s-2,\delta+2}$ 
we need $-\frac{3}{2}< 
\delta<-\frac{1}{2}$.  Combine it with the conditions of Proposition 
\ref{prop:power}; we have that
\begin{math}
 \frac{1}{\beta-2}-\frac{3}{2}\leq \delta<-\frac{1}{2}
\end{math} or $\beta-2>1$. Since $\beta=\frac{2}{\gamma-1}$, it results 
in the condition $\gamma<\frac{5}{3}$, which coincides with the condition of 
the existence and uniqueness theorem \cite[Theorem 4.3]{BK9}.

\subsection{The cosmological model.}
\label{sec:cosmological-context}
The initial value problem for the Euler--Poisson system considered so
far concerned the case of an isolated system where, by definition, the
density, as well as the gravitational potential, vanish at spatial
infinity. 
Here, in this section, we start with homogeneous, isotropic solutions, which have
a spatially constant, non-zero mass density and which describe the
mass distribution in a Newtonian model of the universe.
These homogeneous states can be constructed explicitly, and we
consider deviations from such homogeneous states, which then satisfy a
modified version of the Euler--Poisson system. 
We prove the well-posedness of this setting for initial data which
represent spatially periodic deviations from homogeneous states.
Spatially homogeneous, cosmological solutions for the Euler--Poisson
system can be constructed as follows:  we set
\begin{equation*} 
 \rho (t,x):=  \rho(t)    \quad \mbox{and} \quad  v^k (t,x):=
  \frac{\dot R(t) }{R(t)}x^{k}, \quad k=1,2,3,
\end{equation*}
where $R(t)$ is a
positive scalar function to be determined later. 
Then the continuity equation \eqref{eq:EP:0} results in
\begin{equation*}
  \frac{d}{dt} ( R^3\rho )=0.
\end{equation*}
We obtain the homogeneous mass density
\begin{equation*} 
\widehat \rho(t) = R^{-3}(t)C,\quad t\geq 0.
\end{equation*}
The Euler equations result in 
\begin{equation*}
\frac{\ddot R}{R}x^k=-\partial^k\phi,
\end{equation*}
and it remains to determine the function $R(t)$; the scalar function
$\phi$ is a solution to the Poisson equation.
A short computation, for which we use $C=1$ to simplify the calculations, shows that 
$R(t)$ must be a solution of the differential equation
\begin{equation*} 
  \ddot R + \frac{4 \pi}{3} R^{-2} =0 ,
\end{equation*}
which is the equation of radial motion in the gravitational field of a
point mass. 
For a discussion of its solution, we refer the interested reader to
\cite{Rindler_78}.

In order to study the well-posedness of such a homogeneous model
we investigate the time evolution of small deviations from it,
that is, we consider solutions of the Euler--Poisson system
\eqref{eq:EP:0}--\eqref{eq:sec3-applications:6} of the form
\begin{equation*}
v^a=\widehat v^{a} + V^{a},\quad \rho =\widehat \rho + \sigma ,\quad \phi = \widehat \phi + \Phi .
\end{equation*}
The pressure $p$ is connected to the density by an unspecified,
smooth, equation of state
$p = f( \widehat \rho+ \sigma)-f(\widehat \rho)=: g(\rho,t)$.
We require 
\begin{equation}
\label{eq:sec3-applications:4}
g'>0, \quad \mbox{and} \quad \widehat \rho + \sigma > 0
\end{equation}
where $g'$ denotes $\frac{dg}{d\rho}$. 
To investigate the perturbations $\Phi,\ \sigma,\ V^a$, which we want to assume
spatially periodic, it is useful to perform the following
transformation of variables:
\begin{equation*} 
\tilde x^k  =  R^{-1} (t) \,x^{k}
\end{equation*}
More details on this transformation are given in
\cite{rein94:_local_vlasov}. Its
necessity stems from the fact that in the original variables the
equations for the deviation $g$ would become explicitly $x$-dependent,
which would exclude the possibility of studying spatially periodic
deviations, a class which seems physically reasonable and is
convenient for our mathematical analysis.

If we compute the system satisfied by $V^{a}$, $\sigma$, $\Phi$ in these
transformed variables and afterward drop the tildes, we obtain the following version of the
Euler--Poisson system which governs the time evolution of deviations
from the homogeneous state:
\begin{align}
  \label{eq:sec3-applications:1}
  A^0(U;t)\partial_t U + \sum_{k=1}^3 A^k(U;t) \partial_k
  U + B(U;t) U&= G(U;t) \\
  \label{eq:sec3-applications:3}
   \Delta\Phi &= 4\pi  R^2 \sigma,
\end{align}
where
\begin{equation}\label{eq:sec3-applications:17}
\begin{gathered}
   U = 
  \begin{pmatrix}
    \sigma\\
    V^{a}
  \end{pmatrix},
  \;\: A^0 = 
  \begin{pmatrix}
     \frac{1}{(\widehat \rho+\sigma)^2}{\sc g'}&0\\
    0& \delta_{ab}
  \end{pmatrix},
  \;\:
    A^k = 
  \begin{pmatrix}
    \frac{1}{R} \frac{1}{(\widehat \rho+\sigma)^2}{\sc g'}    {\sc
      v^k}&\frac{1}{R}\frac{1}{\widehat \rho +\sigma}{\sc g'} \delta{}^k_{i} \\
    \frac{1}{R}\frac{1}{\widehat \rho+\sigma}{\sc g'} \delta{}^k_{j}& \frac{1}{R}{\sc v^k}\delta_{ij} 
  \end{pmatrix},
\\[0.2cm]
   B = 
  \begin{pmatrix}
    {\sc 3}\frac{R'}{R}\frac{1}{\widehat \rho + \sigma}    {\sc g'}&0\\
    0&\frac{R'}{R} \delta_{ij}
  \end{pmatrix},
  \;\:
  G = 
  \begin{pmatrix}
    0\\
    -\partial^a\Phi
  \end{pmatrix}.
  \end{gathered}
\end{equation}
Here $\delta_{ij}$  denotes the Kronecker delta.
The coefficients in these matrices satisfy the conditions
\eqref{eq:sec3-applications:4} and therefore $A^0$ is positive definite.

  Several remarks are in order. 
  On one hand, the system \eqref{eq:sec3-applications:1} and
  \eqref{eq:sec3-applications:3} is simpler than the one describing
  the isolated body, since we don't need to use the Makino variable. 
  On the other hand, due to the reference solution, $A^0\neq \Id$, the
  system contains a lower order term $B$ which depends on $t$ 
  moreover, $A^{0},A^{k}$ depend on $t$ as well. 
  The presence of the $A^0$ term complicates somehow the energy
  estimates, but not considerably. 
  Finally, we remark that since the equation of state in that context
  is of a more general nature, we do not use fractional Sobolev spaces.
  The periodicity of the initial data is expressed by using Sobolev
  spaces which were defined over the flat torus $\setT^3$, whose
  corresponding norm is given by
\begin{equation*}
  \Vert u \Vert_{H^m(\setT^3)}^{2} =  \sum_{\alpha\leq m}
\quad  \int_{{\setT^3} } 
  | \partial^{\alpha} u |^2  
  dx, \quad m\in \setN;
\end{equation*}
here $\alpha$ denotes the multi index, a notation we will use in what follows.
In this setting local existence and uniqueness has been shown
\cite{ICH2} for periodic initial  conditions for $\sigma_0, V^a_{0}$ under 
the additional condition that $\int_{\setT^3}\sigma dx=0$.

So we are now in a position to present the continuity of the flow map.

\begin{thm}[Continuity of the flow map  of the Euler--Poisson system in a cosmological context]
\label{thr:sec3-applications:3}
  Let $m\!\!>\!\!\frac{5}{2}$ and
  $(\sigma, V^a)\!\!\in\!\! C^{0}([0,T];H^m(\setT^{3}))\cap
  C^1([0,T];H^{m-1}(\setT^{3}))$ be the solution of the cosmological
  Euler--Poisson system
  \eqref{eq:sec3-applications:1}--\eqref{eq:sec3-applications:3} with
  initial data $(\sigma_0, V^a_0)\in H^m(\setT^{3})$. 
  Let
  \begin{equation*}
  \lim_{n\to \infty}\|(\sigma_0^n,(V_0^a)^n)-(\sigma_0, V_0^a)\|_{H^m(\setT^{3})}\to  0.
\end{equation*}
  Then for sufficiently large $n$, the solution $(\sigma^n, (V^a)^n)$ to
  \eqref{eq:sec3-applications:1}--\eqref{eq:sec3-applications:3}
  exists in the interval $[0,T]$ and
  \begin{equation*}
    \lim_{n\to\infty} \sup_{0\leq t\leq
      T}\|(\sigma^n,(V^a)^n)-(\sigma,V^a)\|_{H^m(\setT^{3})}\to 0.
\end{equation*}
\end{thm}

\begin{proof}
  The proof has some similarities to that of Theorem
  \ref{thr:sec3-applications:2}, but the evolution system in question
  is slightly more complicated than those in Sections
  \ref{sec:case-comp-supp} and \ref{sec:case-density-which} since it
  contains a lower order term $B$ which depends on $t$ and since it
  contains $A^0\neq \Id$. 
  Moreover, $A^{0},A^{k}$ depend on $t$ as well. 
  This will force us to generalize  Theorem \ref{thm:1} to the system
\begin{equation*}
  \label{eq:sec3-applications:15}
  \tag{QSA0}
\begin{cases}
           A^0(U;t)\partial_t U+ \sum_{a=1}^3 A^{a}(U;t)\partial_a U 
           =\widehat G(U;t),\\
           U(0,x)=u_0(x),
         \end{cases}
     \end{equation*}
where $U$ is a function on the torus $\setT^3$ and 
$A^0$ is a uniform positive definite matrix satisfying  
\begin{equation}
\label{eq:sec3-applications:13}
C^{-1} |v|^2\leq v\cdot A^0(u;t)v \leq C |v|^2
\end{equation}
for a positive constant $C$, any vector $v\in \setR^N$ and for all
$ t \in [0,T]$. 
In order to simplify the notation, the function
$\widehat G(U;t)=G(U;t)-B(U;t)U$ from equation
\eqref{eq:sec3-applications:17}.

We then observe that the right-hand side term $\widehat G$ is a sum of
a smooth function $B(U;t)U$ and the source term
$G(U)=-(0,\nabla \Phi)$, where $\Delta \Phi=\sigma$. 
Since $B(U;t)U$ is a smooth function, it then follows from the
difference estimate \eqref{eq:tools:6} that $D_x(B(U;t)U$ is Lipschitz in
$H^{m-1}(\setT^3)$. 
Since we assume that $\int_{\setT^3}\sigma dx =0$, then the Poisson
equation $\Delta \Phi =\sigma$ has a solution for $\sigma \in H^m(\setT^3)$. 
Moreover,
\begin{equation*}
 D_x \nabla \Delta^{-1}: H^{m-1}(\setT^3)\cap 
\bigg\{\int_{\setT^3}\sigma dx=0\bigg\}\to H^{m-1}(\setT^3)\cap 
\bigg\{\int_{\setT^3}\sigma dx=0\bigg\}
\end{equation*}
is a bounded linear operator, and hence Lipschitz.
So we conclude that the right-hand side of \eqref{eq:sec3-applications:15} 
satisfies the condition 
\begin{equation*}
 \|D_x(\widehat G(U)-\widehat G(\widetilde 
U))\|_{H^{m-1}{(\setT^3)}}\leq L
\|U-\widetilde U\|_{H^{m-1}{(\setT^3)}}.
\end{equation*}

It remains to present and prove a  generalization of Theorem \ref{thm:1} 
to the system \eqref{eq:sec3-applications:15}, namely
Theorem \ref{thr:sec3-applications:4}, which will be presented below.
\end{proof}

\subsubsection{The continuity of the flow map for a quasilinear symmetric
  hyperbolic system $A^0\ne \Id$.}
\label{sec:proof-theorem-}

\begin{thm}[The continuity of the flow map for a quasilinear symmetric
  hyperbolic system $A^0\ne \Id$]
\label{thr:sec3-applications:4}
  Let $m>\frac{3}{2}+1$ and assume that
  \eqref{eq:assump:1}--\eqref{eq:assump:3} hold. 
  Let $u_0\in H^m(\setT^3)$ and let
  $U(t)\in C^0([0,T];H^m(\setT^3))\cap C^1([0,T];H^{m-1}(\setT^{3}))$ be the corresponding
  solution to \eqref{eq:sec3-applications:9} with initial data $u_0$. 
  If $\|u_0^n -u_0\|_{H^{m}(\setT^{3})}\to 0$, then for large $n$ the solutions
  $U^n(t)$ to \eqref{eq:sec3-applications:9} with initial data $u_0^n$ exist for
  $t\in[0,T]$, moreover,
  \begin{equation}
    \label{eq:sec3-applications:11}
    \lim_{n\to\infty}\sup_{0\leq t\leq T} \Vert U^n(t)-U(t) \Vert_{H^{m}(\setT^{3})}=0.
  \end{equation}
\end{thm}

The proof of  Theorem  \ref{thr:sec3-applications:4} involves generalizing the low regularity
energy estimate (Lemma \ref{lem:sec1-outline:2}), the continuous
dependence of the solution on the coefficients (Lemma~\ref{lem:sec1-outline:4}), and finally the splitting argument 
outlined below Remark \ref{rem:Katos-example} in Subsection
\ref{sec:continuity-flow-map}.

So the low regularity estimate has the following form.
\begin{lem}[Low regularity energy estimate with $A^{0}\neq \Id$]
\label{lem:sec3a-proof-lemma1A:1}
  Let $m>\frac{3}{2}+1$,\linebreak $A^0$ satisfies
  \eqref{eq:sec3-applications:13}, and
  $A^0,A^a \in L^\infty([0,T]; H^{m}(\setT^3)) $, $\partial_t A^0 \in 
L^\infty([0,T]; L^\infty(\setT^3)) $,
  $F \!\!\!\in\! \!L^\infty([0,T]; H^{m-1}(\setT^3)) $, and $u_0\!\!\in\! H^{m-1}(\setT^3)$.
  Assume $U(t)\!\!\in\! L^{\infty}([0,T] ;H^{m-1}(\setT^3))$ is a solution 
to
  the initial value problem
  \begin{equation*}
    \label{eq:sec3-applications:9}
    \tag{LSA0}
      \begin{cases}
        A^0(t,x)\partial_t U+ \sum_{a=1}^3 A^{a}(t,x)\partial_a U 
=F(t,x),\\
        U(0,x)=u_0(x).
      \end{cases}
  \end{equation*}
  Then for $t\in [0,T]$,
  \begin{equation}
    \label{eq:sec3-applications:14}
    \Vert U(t) \Vert_{H^{m-1}(\setT^3)}^2\leq e^{\int_0^t 
a_m(\tau)d\tau} 
    \bigg( \Vert u_0 \Vert_{H^{m-1}(\setT^3)}^2+ \int_0^t     
\Vert   F(\tau,\cdot)\Vert_{H^{m-1}(\setT^3)}^2 d\tau \bigg),
  \end{equation}
  where
  $a_m(\tau):= C\|A^0(\tau)\|_{H^{m}} \sum_{a=1}^d\Vert
    A^a(\tau)\Vert_{H^{m}}+\|\partial_t A^0\|_{L^\infty}$.  
\end{lem}

\begin{proof}[Proof of Lemma \ref{lem:sec3a-proof-lemma1A:1}]
  We wish to emphasize that we face the same regularity problem
  we described already in the section following Remark
  \ref{rem:sec1-cont-flow-hs:1}, that is, we first have to assume one
  degree more of regularity, and then do an appropriate approximation
  argument. 
  So we would need first to prove an equivalent to Proposition~\ref{lem:sec1:2}.
  However, since we have done this in great detail in Subsection
  \ref{sec:useful-lemmas}, we will just outline the main steps.
  Besides this we proceed in a manner very similar to in the proof
  of Lemma \ref{lem:sec1-outline:2}, especially the idea of using a modified
  commutator estimate.
  However, the presence of the $A^0$ requires the following modification
  (\cite{majda84:_compr_fluid_flow_system_conser}) of the inner product.
  We set
  \begin{equation}
    E_{m}(t)=\sum_{|\alpha|\leq m}^{} \langle 
      \partial^{\alpha}U,A^0\partial^{\alpha}U \rangle_{L^{2}}.
  \end{equation}
Note that 
\begin{equation}
\label{eq:sec3a-proof-lemma1A:17}
\frac{1}{C} \Vert U \Vert_{H^m(\setT^3)}^2 \leq E_m(t) \leq C 
\Vert U \Vert_{H^m(\setT^3)},
\end{equation}
where $C$ is the constant of inequality \eqref{eq:sec3-applications:13}. In 
order to shorten the notation, we shall write $\Vert \cdot \Vert_{H^{m}} $
instead of $\Vert \cdot \Vert_{H^m ( \setT^3 )}$.

Since we consider solutions in $H^{m-1}$, we shall estimate 
$E_{m-1}(t)$. So as in the proof of Proposition \ref{lem:sec1:2}, we 
differentiate the energy and obtain
\begin{align}
  \label{eq:sec3a-proof-lemma1A:2}
  \frac{1}{2}\frac{d}{dt}E_{m-1}(t)
  &= \sum_{|\alpha|\leq m-1}^{}
    \langle
    \partial^{\alpha}U,A^0\partial_{t}\partial^{\alpha}U
    \rangle_{L^{2}}
    +\frac{1}{2} \sum_{|\alpha|\leq m-1}^{} \langle 
\partial^{\alpha}U,\partial_{t}(A^0  )\partial^{\alpha}U 
\rangle_{L^{2}}.
\end{align}

Since  
$\partial_tA^0\in L^{\infty}( [ 0,T ];L^\infty)$,
we obtain 
\begin{equation}
  \label{eq:sec3a-proof-lemma1A:14}
  \sum_{|\alpha|\leq m-1}^{} 
  \langle \partial^{\alpha}U,\partial_{t}(A^0  
)\partial^{\alpha}U \rangle_{L^{2}}
  \leq C \Vert \partial_tA^0 \Vert_{L^{\infty}}\Vert U  \Vert_{H^{m-1}}^{2}.
\end{equation}

The essential part now is to estimate the first term on the right-hand
side of \eqref{eq:sec3a-proof-lemma1A:2}. 
To do this we multiply first the equation
\eqref{eq:sec3-applications:9} by $( A^0 )^{-1}$ and obtain
\begin{equation}
\label{eq:sec3a-proof-lemma1A:3}
 \partial_tU+\sum_{a=1}^3 ( A^0 )^{-1}A^a\partial_a U = ( A^0 )^{-1}F.
\end{equation} 
Then  for each multi-index $\alpha$, $|\alpha|\leq m-1$,  we apply the 
operator 
 $\partial^\alpha$ to \eqref{eq:sec3a-proof-lemma1A:3},  resulting with the 
following equation for  $\partial^\alpha U$:
\begin{equation}
 \partial_t\partial^{\alpha}U+\sum_{a=1}^3 \partial^{\alpha}[ 
(A^0(t,x))^{-1} A^a(t,x)\partial_a U ] =
                       \partial^{\alpha} [( A^0(t,x))^{-1}F(t,x)].
\end{equation} 
We now introduce the identity 
\begin{equation}
\label{eq:sec3a-proof-lemma1A:23}
 (A^0 )^{-1}A^a\partial_a U 
=\partial_a((A^0 )^{-1}A^a U)-\partial_a((A^0 )^{-1}A^a ) U,
\end{equation}
multiply it  by $\partial^\alpha $, and commute the first term on 
the right-hand side of \eqref{eq:sec3a-proof-lemma1A:23} with respect to the 
operator $( \partial^\alpha\partial_a)$, which results in 
\begin{equation}
 \label{eq:sec3a-proof-lemma1A:24}
\begin{split}
 ( \partial^\alpha\partial_a)&[(A^0 )^{-1}A^a 
U]  \\
&=
( \partial^\alpha\partial_a)[(A^0 )^{-1}A^a U]-
(A^0 )^{-1}A^a( \partial^\alpha\partial_a)[ U] + 
(A^0 )^{-1}A^a( \partial^\alpha\partial_a)[U].
\end{split}
\end{equation} 
So now we conclude that for each $\alpha$, $\partial_t\partial^\alpha U$ 
satisfies the equation
\begin{equation}
\label{eq:sec3a-proof-lemma1A:4}
\begin{split}
A^0(t,x)\partial_t[\partial^{\alpha}U]+\sum_{a=1}^3 &
 A^a(t,x)\partial_a[\partial^\alpha U ] 
 \\ &
= -A^0\sum_{a=1}^3 ( \partial^\alpha\partial_a)[(A^0 )^{-1}A^a U]-
(A^0 )^{-1}A^a( \partial^\alpha\partial_a)[ U]\\ 
 &\phantom{=\ }+A^0\sum_{a=1}^3 \partial^\alpha[ \partial_a((A^0 
)^{-1}A^a ) U]+A^0\partial^{\alpha} [( A^0)^{-1}F(t,x)]\\ & =: I +II +III.
\end{split}
\end{equation} 

Now we have to show that each term on the right-hand side of
\eqref{eq:sec3a-proof-lemma1A:4} can be estimated by $\|U\|_{H^{m-1}}$ and
$\|F(t,\cdot)\|_{H^{m-1}}$. 
For the first term we use the third Moser inequality, which is
basically the Kato--Ponce commutator estimate \eqref{eq:kato-ponce} for
$m\in \setN$ (see, e.g., Taylor~\cite[\& 13.3]{taylor97c}). 
The advantage of the identity \eqref{eq:sec3a-proof-lemma1A:24} is
that we use it for the operator $\partial^\alpha\partial_a$, which is of order
$|\alpha|+1$. 
That results in
 \begin{align*}
\label{sec:low-regul-energy-1}
  \|I\|_{L^2}&=     
\|A^0\{(\partial^{\alpha}\partial_a)[A^a U]-A^a
          (\partial^{\alpha}\partial_a)[ 
U]\}\|_{L^2}\\
&  \lesssim \|A^0\|_{L^\infty}
       \{\|D A^a\|_{L^\infty}\|       
U\|_{H^{m-1}}+\|A^a\|_{H^{m}}\|
          U\|_{L^\infty}\}.
  \end{align*}
We now apply the Sobolev embedding theorem to $\|U\|_{L^\infty}$, which allows
  us to conclude that $\|U\|_{L^\infty}\lesssim \|U\|_{H^{m-1}}$, for
  $m-1>\frac{d}{2}$, and similarly for the matrices, so we conclude that
\begin{equation}
 \label{eq:sec3a-proof-lemma1A:25}
 \|I\|_{L^2}\lesssim \|A^0\|_{H^m} \|A^a\|_{H^m}\|U\|_{H^{m-1}}.
\end{equation}   
We turn now to estimate the terms $\|II\|_{L^2}$ and $\|III\|_{L^2}$. 
The inverse matrix $(A^0)^{-1}$ is involved in both of them. 
It follows from assumption \eqref{eq:sec3-applications:13} and from
the Moser estimates \cite[\& 13.3]{taylor97c} that
$\|(A^0)^{-1}\|_{H^m}\leq C_m \|A^0\|_{H^m}$. 
Hence by the multiplication property \eqref{eq:mult}, we obtain 
\begin{equation*}
\begin{split}
 \|II\|_{L^2}= \bigg\| A^0\sum_{a=1}^3 \partial^\alpha[ 
\partial_a((A^0)^{-1}A^a ) U]\bigg\|_{L^2}&\lesssim
\|A^0\|_{L^\infty}\bigg\|\sum_{a=1}^3\partial_a((A^0 
)^{-1}A^a ) U\bigg\|_{H^{m-1}}, \\
&\lesssim
\|A^0\|_{L^\infty}\sum_{a=1}^3\|\partial_a((A^0 
)^{-1}A^a )\|_{H^{m-1}} \| U\|_{H^{m-1}}\\
&\lesssim \| A^0\|_{H^m}^2\bigg(\sum_{a=1}^3\|A^{a}\|_{H^m}\bigg)
\| U\|_{H^{m-1}}.
\end{split}
\end{equation*} 

Similarly, 
\begin{equation*}
 \| III \|_{L^2}=\|( A^0)^{-1}F(t,\cdot)\|_{L^2}\lesssim 
\|A^0\|_{H^m}^2\|F(t,\cdot)\|_{H^{m-1}}.
\end{equation*} 
 
Since the matrices $A^a$ are symmetric, we can apply a  standard integration by 
parts argument to show that
\begin{equation}
 \label{eq:sec3a-proof-lemma1A:28}
 |\langle \partial^\alpha U,A^a\partial_a(\partial^\alpha 
U)\rangle|\lesssim \|\partial_a A^a\|_{L^\infty}\|\partial^\alpha 
U\|_{L^2}^2\lesssim \| A^a\|_{H^m}\|U\|_{H^{m-1}}^2.
\end{equation} 

So now by inequalities
\eqref{eq:sec3a-proof-lemma1A:25}--\eqref{eq:sec3a-proof-lemma1A:28}
and the Cauchy--Schwarz inequality, we obtain that for each $|\alpha|\leq m-1$
\begin{equation}
 \label{eq:sec3a-proof-lemma1A:29}
 \begin{split}
 | \langle
    \partial^{\alpha}U,A^0\partial_{t}\partial^{\alpha}U
    \rangle_{L^{2}}| & \lesssim \sum_{a=1}^3 
    | \langle
    \partial^{\alpha}U,A^a\partial_{a}(\partial^{\alpha}U) 
\rangle_{L^2}|
    + |\langle 
 \partial^\alpha U, I+II +III\rangle_{L^2} | \\  
 & \lesssim \|A^0\|_{H^m}\bigg(\bigg(\sum_{a=1}^3 \|A^a\|_{H^m}\bigg)\|U\|_{H^{m-1}}^2
 + \|U\|_{H^{m-1}}\|F(t,\cdot)\|_{H^{m-1}}\bigg) \\  
& \lesssim \|A^0\|_{H^m}\bigg(\bigg(\sum_{a=1}^3 
\|A^a\|_{H^{m}} +1\bigg)\|U\|_{H^{m}}^2 + 
\|F(t,\cdot)\|_{H^{m-1}}^2\bigg). 
 \end{split}
\end{equation} 

Using the equivalence of norms \eqref{eq:sec3a-proof-lemma1A:17},  equation 
\eqref{eq:sec3a-proof-lemma1A:2}, and inequalities 
\eqref{eq:sec3a-proof-lemma1A:14} and \eqref{eq:sec3a-proof-lemma1A:29}, we obtain
\begin{equation}
 \frac{1}{2}\frac{d}{dt}E_{m-1}(t)\leq C E_{m-1}(t)+ \|F(t,\cdot)\|_{H^{m-1}}^2,
\end{equation} 
where $C$ depends on $\|A^0\|_{H^m}$, $\|A^a\|_{H^m}$ and
$\|\partial_tA^0\|_{L^\infty}$. 
Hence \eqref{eq:sec3-applications:14} holds by the Gronwall inequality
for a solution with one more degree of regularity and the rest proceeds as in
the proof of Lemma \ref{lem:sec3a-proof-lemma1A:1}.

\end{proof}

We turn now to the splitting argument.

The main idea of the proof is to show convergence in the
$H^{m-1}(\setT^3)$ norm for the derivatives of the solution. 
So we start by differentiating equation \eqref{eq:sec3-applications:15} with respect to the
spatial variables $D_x$ and obtain
\begin{equation}
    \begin{cases}
    \begin{split}
       D_x(A^0(U;t))\partial_t U+ A^0(U;t)&\partial_t
      (D_xU)+\textstyle\sum_{a=1}^3A^a(U;t)\partial_a(D_xU)\\
      &+\textstyle\sum_{a=1}^3
      D_x(A^a(U;t))\partial_a U
       = D_x(\widehat  G(u;t)) 
       \end{split}\\
       D_xU(0,x)=D_xu_0(x).
    \end{cases}
\end{equation}

Following the same procedure that was presented in Subsection
\ref{sec:continuity-flow-map}, we obtain that $D_xU$ satisfies
\begin{equation}
  \label{eq:split:1}
  \begin{split}
    A^0(U;t)\partial_t&(D_xU) +\sum_{a=1}^3A^a(U;t)\partial_a(D_xU)\\
    &= D_x(\widehat G(U;t))
    -\sum_{a=1}^3 D_x(A^a(U;t))\partial_a U \\
    &\phantom{=\ }+
    D_x(A^0(U;t))(A^0(U;t))^{-1}\bigg[ \sum_{a=1}^3
      A^a(U;t)\partial_a U -\widehat G(U;t)\bigg].
  \end{split}
\end{equation}

We now set 
\begin{align*}
  H^n &= -\sum_{a=1}^3
       D_x(A^a(U^n;t))\partial_a U^n \\
       &\phantom{=\ }+
       D_x(A^0(U^n;t))(A^0(U^n;t))^{-1}\bigg[
       \sum_{a=1}^3 A^a(U^n;t)\partial_a U^n -\widehat G(U^n;t)\bigg],\\
  H &= -\sum_{a=1}^3
      D_x(A^a(U;t))\partial_a U \\&\phantom{=\ }+
      D_x(A^0(U;t))(A^0(U;t))^{-1}\bigg[
      \sum_{a=1}^3 A^a(U;t)\partial_a U -\widehat G(U;t)\bigg],
\end{align*}
Next, we split $D_{x}U^n=W^n+Z^n$  analogously to Step 2 in the proof of
Theorem \ref{thm:1}. 
Here $W^n$ satisfies the initial value problem
\begin{equation}
  \label{eq:split:3}
    \begin{aligned}
      A^0(U^n;t)\partial_t W^n + \sum_{a=1}^3 A^a (U^n;t) \partial_a W^{n} = H
      +D_x\widehat G(U;t),\\
       W^n(0,x)=D_xu_{0}(x),
    \end{aligned}        
\end{equation}
while  $Z^n $ satisfies
\begin{equation}
  \label{eq:split:4}
    \begin{aligned}
      & A^0(U;t)\partial_t Z^n + \sum_{a=1}^3 A^a (U^n;t) \partial_a Z^{n} = H^n-H +
      D_x(\widehat G(U^n))- D_x(\widehat G(U)),\\
      & Z^n(0,x)=D_xu_{0}^n(x)-D_xu_0(x).
    \end{aligned}        
\end{equation}

The next step consists in proving $W^n\to D_xU$ in the $H^{m-1}$ norm. 
For this, we need a generalization of Lemma \ref{lem:sec1-outline:4} to
the case $A^0\neq\Id$ and which we then have to apply to the system
\eqref{eq:split:3}.
For the sake of brevity we have neither stated this generalization nor
presented its proof, since it is basically a repetition of the same
arguments and applying the low regularity energy estimate.
Moreover, the steps to prove
\begin{equation*}
\lim_{n\to\infty}\sup_{0\leq t\leq T}  \|W^n-D_{x}U\|_{H^{m-1}(\setT^3)}=0,
\end{equation*}
that is, equations
\eqref{eq:sec1-cont-flow-hs:10}, \eqref{eq:sec1-cont-flow-hs:8}, and
\eqref{eq:sec2-main-result:4}, are so similar to those in the proof of
Theorem \ref{thm:1}, that we omit them also.

It remains to show that $Z^n\to 0$ in the $H^{m-1}(\setT^3)$ norm. 
For that, we have to perform steps similar to those which lead to
\eqref{eq:proof:13}, especially Proposition \ref{lem:sec4-proof:1},
which again we will not prove in detail but emphasize apparent
difficulties.
Note that the terms $H^n$ and $H$ are similar to their corresponding
terms in equation \eqref{eq:sec4-proof:5}, but here we have an
additional term, namely, the one in the second row of equation~\eqref{eq:split:1}. 
This additional term is a sum of smooth functions and a function in
the form $F(u)\widehat G(u)$, where $\widehat G(u)$ is the non-linear source term and
$F(u)$ is a smooth function of $u$.
It is not difficult to show that if $D\widehat G(u)$ is Lipschitz in the
$H^{m-1}(\setT^3)$ norm, then $F(u)\widehat G(u)$ is also Lipschitz in the same norm. 
Therefore, in a way similar to the proof of Proposition
\ref{lem:sec4-proof:1} we obtain 
\begin{equation*}
  \|H^n-H\|_{H^{m-1}(\setT^3)}\lesssim 
  \{\|U^n-U\|_{H^{m-1}(\setT^3)}+\|W^n-DU\|_{H^{m-1}
      (\setT^3)}+
    \|Z^n \|_{H^{m-1}(\setT^3)}\}.
\end{equation*}
Finally, applying the low regularity estimate
\eqref{eq:sec3-applications:14} to the initial value problem
\eqref{eq:split:4}, we obtain 
\begin{equation*}
 \lim_{n\to\infty}\sup_{0\leq t\leq T^{\ast}} 
\|Z^n(t)\|_{H^{m-1}(\setT^3)}=0,
\end{equation*}
which completes the proof.  

\biblio

\appendix

\section{Mathematical tools}
\label{sec:mathematical-tools}
Here we list commonly known mathematical tools that are needed for our proofs. 

One of the basic tools for obtaining energy estimates is the
Kato--Ponce commutator estimate. 
Here we shall use the following pseudodifferential operator version
of it (see \cite[\S3.6]{Taylor91}): 
\begin{prop}[The Kato--Ponce commutator estimate]
  \label{lem:sec1-math-prel:3}
  Let $P$ be a pseudo\-differential operator in the class $OPS^s_{1,0}$.
  Then
  \begin{equation}
    \label{eq:kato-ponce}
    \|P(fg)-fP(g)\|_{L^2}\leq C\{\|D
      f\|_{L^\infty}\|g\|_{H^{s-1}}+\|f\|_{H^s}\|g\|_{L^\infty}\},
  \end{equation}
  for any $f\in H^s\cap C^1$ and $g\in H^{s-1}\cap L^\infty$.
\end{prop}

\begin{prop}[Gronwall's inequality]
  \label{lem:sec1-math-prel:2}
  Let $a(t)$ and $b(t)$ be nonnegative and integrable functions in
  $[0,T]$. 
  Suppose $y(t)$ is nonnegative, continuous and satisfies the
  differential inequality
  \begin{equation*}
    \label{eq:sec1-outline:52}
    \frac{d}{dt} y(t) \leq a(t) y(t) + b(t)
  \end{equation*}
  for $t\in[0,T]$. 
  Then
  \begin{equation}
    \label{eq:sec1-outline:53}
    y(t) \leq e^{\int_{0}^{t}a(\tau) d\tau} \bigg( y(0) + \int\limits_0^t     
      b(\tau)d\tau \bigg).
  \end{equation}
\end{prop}

\subsection{Calculus in the  $H^s$ spaces.}
\label{sec:tools-hs-spaces}

\begin{prop}
  \label{prop:H^s}\

  \begin{compactenum}[\upshape (a)]
    \item \label{item:sec4-appendix-H:10} \textbf{Interpolation theorem}
    (\cite[Proposition 1.52]{Bahouri_2011}):
    
    Let $0<s'<s$. Then
    \begin{equation}
      \label{eq:tools:5}
      \|u\|_{H^{s'}}\leq \|u\|_{H^{s}}^{\frac{s'}{s}} 
      \|u\|_{L^2}^{1-\frac{s'}{s}}.
    \end{equation}

    \item \label{item:Multiplication} \textbf{Multiplication}(\cite[Lemma
    A.8]{tao06:_nonlin}):

    Let $s\leq \min\{ s_1,s_2\}$, $s+\frac{d}{2}<s_1+s_2$, $0\leq s_1+s_2$. 
    If $u\in H^{s_1}$ and $v\in H^{s_2}$, then
    \begin{equation}
      \label{eq:mult}
      \|uv\|_{H^{s}}\lesssim 
      \|u\|_{H^{s_1}}
      \|v\|_{H^{s_2}}.
    \end{equation}

    \item \label{item:sec4-appendix-H:11} \textbf{Nonlinear estimates}(\cite[Lemma
    A.9]{tao06:_nonlin}):

    Let $u\in H^s\cap L^{\infty}$, $s\geq0$, and assume $u(x) \in \bar{B}\subset \setR^m$. 
    Let $k$ be the smallest integer greater than $s$ and let
    $F:\setR^m\to\setR^l$ be a $C^{k}$-function such that $F(0)=0$. Then
    $F(u)\in H^s$ with a bound of the form
    \begin{equation}
      \label{nonlinear}
      \|F(u)\|_{H^{s}}\leq 
      C(\|F\|_{C^k(\bar{B})},\Vert u \Vert_{L^{\infty}} )
      \|u\|_{H^{s}}. 
    \end{equation}

    \item
    \label{prop:diff-estimate}
    \textbf{Difference estimate:} Let $u,v\in H^s\cap L^{\infty}$,
    $s\geq0$ and assume \linebreak$u(x),v(x)\in \bar{B}\subset \setR^m$. 
    Let $k$ be the smallest integer greater than $s+1$ and let
    $F:\setR^m\to\setR^l$ be a $C^{k}$-function. Then
    \begin{equation}
      \label{eq:tools:6}
      \|F(u)-F(v)\|_{H^s} \leq C (1+\|u\|_{H^s} 
      +\|v\|_{H^s}) \|u-v\|_{H^s},
    \end{equation}
    where $C=C(\|D_uF\|_{C^{k-1}(\bar{B})},\|u\|_{L^\infty},\|v\|_{L^\infty})$.
    
\textbf{Sketch of the proof:} 
    We write
    \begin{equation*}
      \label{eq:sec4-appendix-H:1}
      F(u)-F(v)= K(u,v)(u-v),
    \end{equation*}
    where
    \begin{equation*}
      \label{eq:sec4-appendix-H:4}
      K(u,v)= \int_0^1 D_u F( \tau u + ( 1-\tau )v  )d\tau.
    \end{equation*}
    Then by \eqref{eq:mult},
    \begin{equation*}
      \|F(u)-F(v)\|_{H^{s}}\lesssim 
      \|K(u,v)\|_{H^{s}}\|u-v\|_{H^{s}},
    \end{equation*}
    and by the above nonlinear estimate we obtain
    \begin{equation*}
    \|K(u,v)\|_{H^{s}}\leq
    C(\|D_uF\|_{C^{k-1}(\bar{B})},\|u\|_{L^\infty},\|v\|_{L^\infty})
    (\|u\|_{H^{s}}+\|v\|_{H^{s}}).
  \end{equation*}
 
    \item \textbf{A power estimate:} The following was proved by Kateb
    \cite{kateb03:_besov} (see also Runst and Sickel
    \cite{runst96:_sobol_spaces_fract_order_nemyt}).
    Let $u\in H^{s}\cap L^\infty$, $1<\beta$, $ 0<s<\beta +\frac{1}{2}$. Then
    \begin{equation}
      \label{eq:tools:7}
      \||u|^\beta\|_{H^{s}}\leq  C(\|u\|_{L^\infty})   \|u\|_{H^{s}}.
    \end{equation}
  \end{compactenum}
\end{prop}

\subsection{Properties of the weighted fractional Sobolev spaces $H_{s,\delta}$.}
\label{sec:prop-weight-sobol}

For the proofs of the calculus type properties see \cite{BK8}.   
\begin{prop}
  \label{prop:w}\
  \begin{compactenum}[\upshape (a)]

    \item 
    \label{item:sec4-appendix-H:4} \textbf{Embedding}:
    \begin{equation}
      \label{eq:appendix-w:1}
      \|\partial_i u\|_{H_{s-1,\delta+1}}\leq \| u\|_{H_{s,\delta}}.
    \end{equation}

    \item
    \label{item:sec4-appendix-H:5}
    \textbf{Multiplication}: Let $s\leq \min\{ s_1,s_2\}$,
    $s+\frac{d}{2}<s_1+s_2$, $0\leq s_1+s_2$,\linebreak
    $\delta-\frac{d}{2}\leq \delta_1+\delta_2$,
    $u\in H_{s_1,\delta_1}$, and $v\in H_{s_2,\delta_2}$. Then
    \begin{equation}
      \label{eq:mult:weight}
      \|uv\|_{H_{s,\delta}}\lesssim \|u\|_{H_{s_1,\delta_1}}
      \|v\|_{H_{s_2,\delta_2}}.
    \end{equation}
  
    \item \textbf{Nonlinear estimate in the weighted spaces}: Let
    $u\in H_{s,\delta}\cap L^{\infty}$, $s\geq0$, $\delta\in\setR$, and assume $u(x)\in \bar{B}\subset\setR^m$. 
    Let $k$ be the smallest integer greater than $s$ and let
    $F:\setR^m\to\setR^l$ be a $C^{k}$-function such that $F(0)=0$. Then
    $F(u)\in H_{s,\delta}$ with a bound of the form
    \begin{equation}
      \label{eq:14}
      \|F(u)\|_{H_{s,\delta}}\leq 
      C(\|F\|_{C^k(\bar{B})},\Vert u \Vert_{L^{\infty}} )
      \|u\|_{H_{s,\delta}}.
    \end{equation}

    \item \label{item:sec4-appendix-H:6} \textbf{Difference estimate}:
    Let $u,v\in H_{s,\delta}\cap L^{\infty}$, $s\geq0$,
    $\delta\in\setR$, and assume $u(x),v(x)\in \bar{B}\subset \setR^m$. 
    Let $k$ be the smallest integer greater than $s+1$ and let
    $F:\setR^m\to\setR^l$ be a $C^{k}$-function. Then
    \begin{equation}
      \label{eq:tools-w:2} 
      \|F(u)-F(v)\|_{H_{s,\delta}} \leq 
      C(1+\|u\|_{H_{s,\delta}}+\|v\|_{H_{s,\delta}})
      \|u-v\|_{H_{s,\delta}},
    \end{equation}
    where $C=C(\|D_uF\|_{C^{k-1}(\bar{B})},\|u\|_{L^\infty},\|v\|_{L^\infty})$.
  
    \item \label{item:sec4-appendix-H:7} \textbf{Power estimate}: Let
    $ \frac{d}{2}<s<\beta +\frac{1}{2}$. Then
    \begin{equation}
      \label{eq:kateb}
      \||u|^\beta\|_{H_{s,\delta}}\leq  C\|u\|_{H_{s,\delta}}.
    \end{equation}

    \item \label{item:sec4-appendix-H:8} \textbf{Embedding into
      $L^\infty$}: Let $\frac{d}{2}<s$ and $-\frac{d}{2}\leq \delta$. Then
    \begin{equation}
      \label{eq:tools-em-w}
      \|u\|_{L^\infty}\leq C \|u\|_{H_{s,\delta}}.
    \end{equation}
     
    \item \label{item:sec4-appendix-H:9} \textbf{Interpolation}: Let
    $0<s<s^\prime$. Then
    \begin{equation}
      \label{eq:tools:inter-w}
      \| u\|_{H_{s,\delta}}\leq \| u\|_{H_{0,\delta}}^{1-\frac{s}{s^\prime}}\| 
      u\|_{H_{s^\prime,\delta}}^{\frac{s}{s^\prime}}.
    \end{equation}
\end{compactenum}
\end{prop}

The following lemma was proved in \cite{BK9}.

\begin{lem}[Improved multiplication
      estimate]
      \label{lem:improved-mult}
      Let $u_i\in H_{s,\delta_i}$ for $i=1,\ldots,m$ and
    $s>\frac{d}{2}$. 
    If
    \begin{math}
      \label{eq:all-proofs:1}
      \delta\leq \delta_1+\cdots+\delta_m+\frac{(m-1)d}{2},
    \end{math}
    then
    \begin{math}
      \label{eq:all-proofs:2}
      u=u_1\cdots u_m\in H_{s,\delta}
    \end{math}
    and
    \begin{equation}
      \label{eq:tools-w:1}
      \|u\|_{H_{s\delta}}\leq       C\prod_{i=1}^m\|u_i\|_{H_{s\delta_i}}.
    \end{equation}
\end{lem}

\subsection{Weak convergence.}
\label{sec:weak-convergence}
\begin{prop}
  \label{prop:weak-conv}
  Let $X$ and $Y$ be two Hilbert spaces such that $X\subset Y$ and the embedding
  is continuous. 
  Suppose $\{x_n\}$ is a sequence in $X$ and such that $x_n$ converges
  weakly to $\tilde{x}_0$ in $X$ and $x_n$ converges to $x_0$ in the norm
  of $Y$. 
  Then
  \begin{equation*}
    \tilde{x}_0=x_0.
  \end{equation*}
\end{prop}

\begin{proof}
  By the assumptions, $f(x_n)\to f(\tilde{x}_0)$ for all $f\in X'$,
  where $X'$ denotes the dual space. 
  Since $X\subset Y$, we have that $Y'\subset X'$, and hence
  $f(x_n)\to f(\tilde{x}_0)$ for all $f\in Y'$. 
  In addition, $\{x_n\}$ converges in the norm of $Y$ to $x_0$, therefore
  it also converges weakly in $Y$, that is, $f(x_n)\to f({x}_0)$ for all
  $f\in Y'$. 
  Thus $f(\tilde{x}_0)=f(x_0)$ for all $f\in Y'$, which implies that
  $\tilde{x}_0=x_0$.  
\end{proof}

\biblio

\bibliographystyle{amsalpha-url} 
\bibliography{bibgraf}
\end{document}